\documentclass[global]{svjour}

\usepackage{amsmath}
\usepackage{latexsym}
\usepackage{epic}
\usepackage{eepic}
\usepackage{mathrsfs}
\usepackage{times}
\usepackage{array}
\usepackage{pstricks}

\input{boxedeps.tex}

\SetEPSFDirectory{}
\HideDisplacementBoxes

\SetRokickiEPSFSpecial  

\newfont{\Ma}{msam10}

\DeclareMathAlphabet{\ams}{U}{msb}{m}{n}
\DeclareMathAlphabet{\goth}{U}{euf}{m}{n}

\def\H{\ams{H}}
\def\C{\ams{C}}
\def\R{\ams{R}}
\def\Z{\ams{Z}}
\def\Q{\ams{Q}}

\def\lcm{\text{lcm}}

\newpsobject{showgrid}{psgrid}{subgriddiv=1,griddots=5,gridlabels=6pt}

\def\SS{\goth{S}}
\def\AA{\goth{A}}
\def\CC{\mathscr C}
\def\FF{\mathscr F}
\def\LL{\mathscr L}
\def\PP{\mathscr P}

\def\top{\text{top}}
\def\vol{\text{vol}}

\def\v{\mathbf{v}}
\def\e{\mathbf{e}}
\def\u{\mathbf{u}}

\def\aa{\alpha}
\def\ss{\sigma}
\def\bb{\beta}
\def\om{\Omega}

\def\i{\text{i}}
\def\j{\text{j}}
\def\k{\text{k}}

\def\ov{\overline}

\def\isom{\text{Isom}\,}

\newcommand{\1}{\boldsymbol{1}}
\newcommand{\2}{\boldsymbol{2}}
\newcommand{\3}{\boldsymbol{3}}
\newcommand{\4}{\boldsymbol{4}}
\newcommand{\5}{\boldsymbol{5}}
\newcommand{\6}{\boldsymbol{6}}
\newcommand{\7}{\boldsymbol{7}}
\newcommand{\8}{\boldsymbol{8}}
\newcommand{\9}{\boldsymbol{9}}
\newcommand{\0}{\boldsymbol{0}}

\newcommand{\HH}{\text{\bf H}}

\newcommand{\sss}{\scriptstyle}

\title{Coxeter Groups and Hyperbolic Manifolds}

\author{Brent Everitt\thanks{\email{bje1@york.ac.uk}}
} 
\institute{Department of Mathematics, University of York, York
YO10 5DD, England 
}

\titlerunning{}
\authorrunning{Brent Everitt 
}

\begin{document}

\maketitle



\begin{abstract}
The rich theory of Coxeter groups is used to provide an algebraic construction of finite
volume hyperbolic 
$n$-manifolds. Combinatorial properties of finite images of these groups can be used to
compute the volumes of the 
resulting manifolds. Three examples, in $4,5$ and $6$-dimensions, are
given, each of very small volume, and in one case of smallest possible
volume.
\end{abstract}


\section{Introduction}\label{intro}

In the last quarter of a century,
$3$-manifold topology has been revolutionised by Thurston
and his school. 
This has generated a huge literature on hyperbolic $3$-manifolds building on the
classical 
$2$-dimensional case. Balanced against this is a relative paucity
of techniques and examples of hyperbolic $n$-manifolds for $n>3$. 
Recent work of Ratcliffe and Tschantz has provided examples of smallest possible volume
when $n=4$ by identifying the faces of polytopes (\cite{Ratcliffe00}, see \cite{Davis85} for
a seminal example).
On the other hand, in a construction more algebraic in nature, they computed 
the covolumes
of the principal congruence subgroups of level $p$ a prime in $\text{PO}_{1,n}\Z$,
giving an infinite family of manifolds, with known volume, in every dimension \cite{Ratcliffe97} .
In a similar vein, \cite{Everitt00} allows one to construct manifolds in arbitrary dimensions by considering
the kernels of representations of Coxeter groups onto finite classical groups.
However, it turns out that the volumes of the manifolds arising from both constructions
are not particularly small (eg: in \cite{Ratcliffe97} the smallest $6$-dimensional
example has Euler characteristic
$|\chi|=44226$). Thus there seem to be no constructions in arbitrary dimensions
that give examples, in low dimensions at least, of very small volume. 

This paper is an attempt to redress this. A construction is given that essentially makes algebraic the 
classical geometric idea of identifying the faces of a polytope.
Algebraic is meant in the following sense: any finite volume hyperbolic 
manifold arises as a quotient $\H^n/\Pi$ by the action of a group $\Pi$ acting properly, freely, 
cofinitely by isometries on $\H^n$. The fundamental group of the manifold is thus 
isomorphic to $\Pi$.
Moreover, by Mostow rigidity \cite{Margulis70,Mostow73,Prasad73}, finite volume hyperbolic $n$-manifolds,
for $n\geq 3$, are completely determined by their fundamental groups, in that two such manifolds 
are isometric precisely when their fundamental groups are isomorphic. 
Thus, to construct $M=\H^n/\Pi$,
we construct $\Pi$.

The geometric notions of a free proper cofinite action can be replaced by the algebraic conditions that
$\Pi$ is a torsion free lattice in the Lie group $\isom\H^n\cong\text{PO}_{1,n}\R$. 
The task can then be spilt into two parts. 
Start with a Coxeter group $\Gamma$ embedded as a lattice in $\isom\H^n$.
Fortunately, there already
exists an extensive literature on such hyperbolic Coxeter groups (see \cite{Vinberg85} and 
the references there). By Selbergs lemma \cite{Selberg60}, such a $\Gamma$ is virtually torsion free
(as indeed is any abstract Coxeter group), so 
the second step is to find the desired torsion free $\Pi$ as a subgroup of finite index.
The machinery necessary for this second step preoccupies most of this paper--it is
combinatorial in nature and uses only the algebraic properties of $\Gamma$. Basic properties of the
resulting manifolds, such as volume and orientability, can then be determined from the structure of $\Gamma$.

The paper is organised as follows: \S 2 contains preliminary material about Coxeter groups;
\S 3 gives the construction and the important role played by the
conjugacy classes
of finite Coxeter groups; these classes are determined in \S 4, where the
results of Roger Carter's lengthy
investigation \cite{Carter72} are summarized in a form amenable to our purposes;
computing volumes
is achieved via the (standard) techniques of \S 5, and 
\S 6 contains three examples, one in each of $4,5$ and $6$-dimensions. These three are either the 
smallest possible ($n=4$)
or would appear to be the smallest\footnote{After this paper was written, John Ratcliffe and 
Steven Tschantz constructed a $5$-dimensional manifold with $\frac{1}{8}$ the volume of our
example \cite{Ratcliffe03}.}
known ($n=5,6$).

It is a pleasure to acknowledge many helpful discussions with Colin
Maclachlan.


\section{Coxeter Groups}

Let $\Gamma$ be a group and $S\subset\Gamma$ finite. The pair $(\Gamma,S)$ is
called a Coxeter system (and $\Gamma$ a Coxeter group) if $\Gamma$ admits a
presentation with generators
$s_{\alpha}\in S$, and relations,
$$
(s_{\alpha}s_{\beta})^{m_{\alpha\beta}}=1,
$$
where $m_{\alpha\beta}\in\Z^+\cup\{\infty\}$ and $m_{\alpha\alpha}=1$
(see \cite{Bourbaki68,Humphreys90} for basic facts 
about Coxeter groups).
It is customary to omit relations
for which $m_{\alpha\beta}=\infty$.
Associated to any Coxeter group is its symbol, with nodes indexed by the
$s_{\alpha}$, where nodes $s_{\alpha}$ and $s_{\beta}$ are joined by an edge
labelled $m_{\alpha\beta}$ if
$m_{\alpha\beta}\geq 4$, an unlabelled edge if 
$m_{\alpha\beta}=3$ and a dotted edge if $m_{\aa\bb}=\infty$. 
It is a standard abuse of notation to denote a Coxeter
group and its symbol by the same letter.

Coxeter groups satisfy a Freiheitssatz, in the sense that for 
any $S'\subseteq S$, the subgroup of $\Gamma$ generated by the
$s_\aa\in S'$ is also a Coxeter group with symbol obtained from $\Gamma$
by removing those nodes not in $S'$ and their incident edges. 
In general, a subsymbol of $\Gamma$ is
some subset of the vertices and their incident edges.
If a symbol $\Gamma$ is disconnected with
connected components $\Gamma_1,\ldots,\Gamma_k$, then 
$\Gamma$ is isomorphic to $\Gamma_1\times\cdots\times\Gamma_k$. 
A group with connected
symbol is called irreducible.

It should come as no surprise that in finding torsion free subgroups of
Coxeter groups the finite, or spherical, Coxeter groups play a central role.
The finite irreducible Coxeter groups are well known to be:
the Weyl groups of 
simple Lie algebras over $\C$, the dihedral groups, the group of symmetries 
of a regular dodecahedron, and the group of symmetries of the regular 
4-dimensional polytope, the 120-cell. 
The Killing-Cartan
notation will be used throughout this paper, in which the Weyl groups fall into three
infinite classical families,
$$
A_n (n\geq 1), B_n (n\geq 2)\text{ and }D_n (n\geq 4),
$$
as well as five exceptional groups,
$G_2, F_4, E_6, E_7\text{ and }E_8$.
The non-Weyl groups are the dihedral $I_2(m)$ for $m\geq 3$, and the dodecahedron
and 120-cell symmetry groups $H_3$ and $H_4$ (note that there are the
exceptional isomorphisms $A_2\cong I_2(3)$, $B_2\cong I_2(4)$ and
$G_2\cong I_2(6)$). 
%
The symbols for the Weyl groups are in 
Tables \ref{roots1}-\ref{roots2}, and the non-Weyl groups in 
Proposition \ref{nonweyl2} and Table \ref{roots3}.
Given an arbitrary (not necessarily connected) Coxeter symbol, it represents a
finite 
Coxeter group if and only if each of its connected components is in
the list above.

A smaller role is reserved for the 
parabolic Coxeter (or affine Weyl)
groups: the irreducible ones are 
$\tilde{A}_n (n\geq 1)$, $\tilde{B}_n (n\geq 3)$, $\tilde{C}_n (n\geq
2)$, $\tilde{D}_n (n\geq 4)$, $\tilde{E}_6$,  $\tilde{E}_7$,
$\tilde{E}_8$, $\tilde{F}_4$,  and $\tilde{G}_2$. 
For their symbols,
consult \cite[\S 2.5]{Humphreys90}. In general, a Coxeter group is
parabolic (or Euclidean) if and only if its symbol has connected
components corresponding to the groups in this list. 

A convex polytope $P$ in $\H^n$ is an intersection,
$$
P=\bigcap_{s\in S}H_s^-,
$$
of closed half spaces $H_s^-$ bounded by hyperplanes $H_s$. We will assume that $S$
is finite (although we could take the collection of hyperplanes to be infinite but locally finite), 
and to avoid degenerate cases,
that $P$ contains a non-empty
open subset of $\H^n$. A Coxeter system $(\Gamma,S)$ is {\em
hyperbolic\/} of dimension $n$ if and only if there is a $P\subset \H^n$ such that 
assigning $s_\aa\in
S$ to the reflection in the hyperplane $H_{s_\aa}$ induces an
isomorphism between $\Gamma$ and the group generated by reflections in
the faces of $P$. Call $\Gamma$ cocompact (resp. cofinite) if $P$ is
compact (resp. of finite volume) in $\H^n$.

Given an abstract Coxeter 
system $(\Gamma,S)$, when can one realise it as a cocompact/cofinite
hyperbolic Coxeter group?
Dealing first with hyperbolicity, a Gram matrix $G(\Gamma)$ for $\Gamma$
has rows and 
columns indexed by $S$; if $m_{\aa\bb}\in\Z^+$, set the $\aa\bb$-th
and $\bb\aa$-th entries of $G(\Gamma)$ to be
$-\cos(\pi/m_{\aa\bb})$, and if $m_{\aa\bb}=\infty$, choose the $\aa\bb$-th
and $\bb\aa$-th entries to be some real $c_{\aa\bb}$ with
$c_{\aa\bb}\leq -1$. Recalling that a symmetric matrix has signature
$(p,q)$ if it has precisely $p$ eigenvalues that are $<0$ and $q$ that
are $>0$,

\begin{theorem}[\cite{Vinberg85}, Theorem 2.1]\label{Vinberg1}
A Coxeter group is $n$-dimensional hyperbolic if there exist
$c_{\aa\bb}$'s such that $G(\Gamma)$ has signature $(1,n)$.
\end{theorem}

Note that in \cite[\S 6.8--6.9]{Humphreys90},
the sufficient condition of Theorem \ref{Vinberg1} (as well as the
condition that the Gram matrix be non-degenerate) is taken as the
definition of hyperbolic Coxeter group, and so a more restricted class
of groups is obtained, indeed, classified. A classification in the 
general case remains elusive.

The configuration of the 
bounding hyperplanes of $P$
can be described as follows. If $m_{\aa\bb}\in\Z^+$ then the
corresponding hyperplanes intersect with dihedral angle
$\pi/m_{\aa\bb}$; if $c_{\aa\bb}=-1$ they are parallel, that is, intersect at
infinity; and if $c_{\aa\bb}<-1$ they are ultraparallel (they do not
intersect in the closure of $\H^n$) with a common perpendicular geodesic
of length $\eta_{\aa\bb}$, where $-\cosh\eta_{\aa\bb}=c_{\aa\bb}$. 
Because of these last two, some authors prefer to embellish the Coxeter
symbol, replacing the dotted edges by thick solid ones if
$c_{\aa\bb}=-1$, or by ones labelled $\eta_{\aa\bb}$.


To determine if $(\Gamma,S)$ is realised by a
cocompact/cofinite $P\subset\H^n$, let $\FF=\FF(\Gamma,S)$ be the collection of
finite subgroups of $\Gamma$ generated by subsets $S^{\prime}\subset S$
(including $S^{\prime}=\emptyset$, which generates the trivial group).
Partially order $\FF$ by inclusion.
Similarly $\overline{\FF}$ is the poset 
obtained by taking both the finite and the parabolic subgroups. The poset of an abstract
combinatorial polytope $P$ is the set of cells of $P$, partially ordered by inclusion.

\begin{theorem}[\cite{Vinberg85}, Proposition 4.2]\label{Vinberg2}
An $n$-dimensional hyperbolic Coxeter group $(\Gamma,S)$ is cocompact 
(resp. cofinite) iff $\FF$ (resp. $\overline{\FF}$)
is isomorphic as a partially ordered set to the poset of some
$n$-dimensional abstract combinatorial polytope.
\end{theorem}

A slight variation yields a very useful $\Gamma$-complex. Let $\Gamma\FF=
\Gamma\FF(\Gamma,S)$ be the set
of all right cosets
of the subgroups in $\FF$.
Partially order
this set by inclusion and let 
$\Sigma$ be the affine complex realising $\Gamma\FF$:
the $k$-cells of $\Sigma$
are the chains
$\ss_0<\ss_1<\cdots<\ss_k$ with $\ss_i\in \Gamma\FF$. 
The group $\Gamma$ acts cellularly on $\Sigma$ by right multiplication.
If $\Gamma$ is hyperbolic
then $\Sigma$
is just the barycentric subdivision of the tessellation of $\H^n$ by congruent copies of
the polytope $P$.
But in fact for any Coxeter system, the complex turns out to be negatively curved:

\begin{theorem}
Let $(\Gamma,S)$ be a Coxeter system and equip $\Sigma$ with the natural piecewise
Euclidean metric 
where each $k$-cell is
isometric to a regular $k$-dimensional Euclidean simplex.
Then this is a $\text{CAT}(0)$ space with the $\Gamma$ action by right multiplication isometric.
\end{theorem}

For basic facts concerning CAT(0) spaces, \cite{Davis98} is recommended,
where this result is Corollary 6.7.5 (see also 
\cite{Moussong88}).

\section{Torsion free subgroups of Coxeter groups}\label{modules}

The construction of torsion free subgroups of a Coxeter group requires
answers to the following two questions: where exactly is the torsion, and given a subgroup,
when does it avoid it?

\begin{theorem}\label{torsion}
Let $(\Gamma,S)$ be a Coxeter system. Then any element of finite order in $\Gamma$ is conjugate to
an element of a finite subgroup generated by some $S'\subseteq S$.
\end{theorem}

This is exercise V.4.2 in Bourbaki \cite{Bourbaki68}, and an algebraic proof 
using root systems appears in \cite{Howlett93}. We provide a geometrical
proof that is presumably well known, but does not appear to be in the literature:

\begin{proof}
Consider the isometric $\Gamma$-action by right multiplication on the complex
$\Sigma$ above. Any torsion
element $\gamma$ in a group of isometries of a $\text{CAT}(0)$ space has a fixed point 
\cite[Corollary II.2.8]{Bridson99},
and since the $\Gamma$-action is cellular, $\gamma$ must fix a $k$-cell 
$\ss_0<\ss_1<\cdots<\ss_k$ with $\ss_i=\langle S_i\rangle g$ for
$g\in\Gamma$ and $S_i\subset S$ with $\langle S_i\rangle$ finite. Thus,
$g\gamma g^{-1}$ fixes the $k$-cell 
$\langle S_0\rangle<\cdots<\langle S_k\rangle$.\hfill\qed
\end{proof}

To address the second question above, it will turn out to be convenient to use the
language of permutation modules rather than subgroups. Let $\Omega$ be a
finite set and $U$ a vector space over $\C$ with basis $\Omega$. An
action of any group $\Gamma$ on $\Omega$ can be extended linearly to an
action on $U$, thus giving $U$ the structure of a (permutation) $\C\Gamma$-module. 
Call the module transitive if it arises via a
transitive $\Gamma$-action on $\Omega$. 

Conjugacy classes of subgroups $\Pi$ of index $|\Omega\,|$ in a group $\Gamma$ are in
one to one correspondence with equivalence classes of transitive representations
$\Gamma\rightarrow\text{Sym}\Omega$. These in turn are in one to one
correspondence with isomorphism classes of transitive permutation
$\C\Gamma$-modules. 

Let $\CC$ be a 
conjugacy class of torsion elements in $\Gamma$. We will
say that the $\C\Gamma$-module $U$ {\em avoids\/} $\CC$ when the corresponding
conjugacy class of subgroups has empty intersection with $\CC$; $U$ is
{\em torsion free\/} when it avoids the conjugacy classes of all torsion elements.
The following is well known and easily proved: 

\begin{proposition}\label{prop1}
Let $U$ be a $\C\Gamma$-module for some group $\Gamma$ with basis $\Omega$. Then $U$ is torsion
free exactly when it avoids the conjugacy classes of elements of prime
order. Moreover, if $\,\CC$ is such a class, and $\gamma\in\CC$ any
element, then $U$ avoids $\CC$ exactly when $\gamma$ fixes no
point of $\Omega$. 
\end{proposition}




So, if we had a list of representatives for the conjugacy classes of
prime order torsion in $\Gamma$,
verifying that a given $U$ is torsion free would become the simple matter of checking that 
no element in the list fixed a point of $\Omega$. By Theorem \ref{torsion}, such a list
could be compiled in particular for $(\Gamma,S)$ a Coxeter group by listing
representatives of the conjugacy 
classes in the finite subgroups
generated by $S'\subseteq S$. Indeed, we can restrict to the finite subgroups generated
by $S'$ that are maximal 
with the property that they generate a finite subgroup.
The list of torsion elements so obtained will very probably include
redundancies, but will certainly be complete, which is clearly all that
matters for the task at hand.
A discussion of conjugacy in finite Coxeter groups is the subject of $\S 4$.

Rather than dealing with an indigestible whole, we can divide the task of finding a 
torsion free $U$ into manageable pieces.
Given $\C\Gamma$-modules $U_i,i=1,2$, with bases $\Omega_i$, let
$\Gamma$ act on $\Omega_1\times\Omega_2$ via
$(\u,\v)^{\gamma}=(\u^{\gamma},\v^{\gamma})$. Extending this action linearly to
the complex vector space with basis $\Omega_1\times\Omega_2$ gives a
permutation $\C\Gamma$-module that may be identified with $U_1\otimes
U_2$. The problem is that $U_1\otimes U_2$ may not be transitive. 
Choose an orbit in the $\Gamma$-action on $\Omega_1\times\Omega_2$,
and extend linearly to a $\Gamma$-action on the $\C$-vector space
with basis this orbit. 
By an abuse of notation, we will use
$U_1\otimes U_2$ to denote any one of the resulting (transitive permutation)
$\C\Gamma$-modules (in the examples in \S 6, $U_1\otimes U_2$ is always transitive).

The following follows by definition and Proposition \ref{prop1}

\begin{lemma}\label{lemma3.1}
Let $U_i, i=1,2$, be $\C\Gamma$-modules and $\CC$ a conjugacy class of
torsion elements in the group $\Gamma$. Then $U_1\otimes U_2$ avoids $\CC$ if and
only if at least one of the $U_i$ does.
\end{lemma}

Thus we may build a torsion free $\C\Gamma$-module $U$ by finding
$U_1,\ldots,U_k$ 
such that each conjugacy class of torsion elements is avoided by at
least one of the $U_i$, and then letting $U=\otimes U_i$. 

In order to determine the
volume of the resulting manifold we will need to be able to show that
$\otimes U_i$ is a transitive module. The remaining results in this
section give us the tools to do this, with the first one just an
elementary result in group theory rewritten in our language,

\begin{lemma}\label{lemma3.2}
Let $U_i, i=1,2$ be transitive $\C\Gamma$-modules with bases $\Omega_i$, and $\Omega$ an orbit 
in $\Omega_1\times\Omega_2$ arising from the $\C\Gamma$-module $U_1\otimes U_2$. Then
$\lcm |\Omega_i|$ divides $|\Omega|$, and $|\Omega|\leq \prod |\Omega_i|$, with equality if
the $|\Omega_i|$ are relatively prime.
\end{lemma}

\begin{proof}
Intersect $\Omega$ with the sets $\{u\}\times\Omega_2$ as $u$ ranges over $\Omega_1$. As the $\Gamma$-action
on $\Omega_1$ is transitive there are $|\Omega_1|$ such sets, each of the same size, thus $|\Omega_1|$
divides $|\Omega|$, and similarly for $\Omega_2$, giving that $\lcm|\Omega_i|$ divides $|\Omega|$. The
upper bound is immediate, as is the final statement, for under the conditions given, $\lcm|\Omega_i|=\prod|\Omega_i|$.
\hfill\qed
\end{proof}

Recall that a $\Gamma$-action on $\om$ is imprimitive if $\om=\bigcup \om_i$,
a disjoint union of blocks with the $1<|\om_i|<|\om|$, 
and $\om_i^{\gamma}=\om_j$ for all $i$ and $\gamma\in\Gamma$
(ie: the action on $\om$ induces an action on the $\om_i$). Call a
module $U$ imprimitive if the $\Gamma$-action on its basis is.

\begin{lemma}\label{lemma3.3}
Let $U$ be an imprimitive $\C\Gamma$-module with basis $\Omega$ and $\overline{\Omega}$ the set
of blocks of imprimitivity with the corresponding $\C\Gamma$-module $\overline{U}$ transitive
(and hence each block having the same size $\mu$).
Let $F$ be a finite subgroup of $\Gamma$ with a fixed point in $\overline{\Omega}$ but no non-identity 
element fixing a point in $\Omega$. Then $|F|$ divides $\mu$.
\end{lemma}

\begin{proof}
If $B$ is the block corresponding to the fixed point in $\overline{\Omega}$ of $F$, then it is a disjoint
union of $F$-orbits, each of which has the form $\{\gamma(x)\,|\,\gamma\in F\}$ for some $x\in B$. 
The condition on the $\Gamma$-action on $\Omega$ gives that each $F$-orbit has size
$|F|$, which thus divides the size of $B$.\hfill\qed
\end{proof}

Notice that we can draw the same conclusion if we merge $\Omega$ and
$\overline{\Omega}$: if $U$ is transitive and $F$ acts without fixed
points on $\Omega$ then $|F|$ divides $|\Omega|$. Consequently, if
$\LL(\Gamma)$ is the lowest common multiple of the orders of the finite
subgroups of $\Gamma$, then $\LL(\Gamma)$ divides the size of the basis
for any torsion free module $U$.
In \S 6 we will use Lemma \ref{lemma3.3} in the following special case,

\begin{proposition}\label{prop2}
Let $U_i,i=1,2$ be transitive $\C\Gamma$-modules with bases $\Omega_i$, $F$ a finite subgroup of
$\Gamma$ with a fixed point in $\Omega_1$ but no non-identity element fixing a point in 
$\Omega_2$. If $\Omega$ is an orbit in $\Omega_1\times\Omega_2$ arising from the $\C\Gamma$-module
$U_1\otimes U_2$, then $|\Omega_1||F|$ divides $|\Omega|$.
\end{proposition}

\begin{proof}
The $\Gamma$-action on $\Omega$ is imprimitive with blocks of imprimitivity the intersections with
$\Omega$ of the sets $\{u\}\times\Omega_2$ as $u$ ranges over $\Omega_1$. Letting $\overline{\Omega}$ be the
set of blocks, the $\Gamma$-action on $\overline{\Omega}$ is transitive as the action on $\Omega_1$ is.
If $F$ fixes $u_0\in\Omega_1$ then it fixes the block $\Omega\cap\{u_0\}\times\Omega_2$ in 
$\overline{\Omega}$, whereas 
if an element of $F$ fixes a point of $\Omega$ then restricting to the second coordinate gives a fixed
point in $\Omega_2$. Thus by Lemma \ref{lemma3.3}, $|F|$ divides the size of each block in $\overline{\Omega}$
of which there are $|\Omega_1|$.\hfill\qed
\end{proof}

Finally, $U$ is {\em orientable\/} whenever there is some vertex that is
not fixed by any word involving an odd number of the generators $S$ (if
$U$ is transitive, then {\em no\/} vertex will be fixed by such a word
as $w$ involves an odd number of generators precisely when 
$w_1ww_1^{-1}$ does for any $w_1$).

\begin{proposition}\label{prop3}
\begin{enumerate}
\item If $\,U_1$ is orientable then so is $U_1\otimes U_2$.
\item If $U$ is transitive and torsion free, with $\Pi$ the
(conjugacy class of the) corresponding subgroup, then the manifold
$\H^n/\Pi$ is orientable if and only if $\,U$ is orientable. 
\end{enumerate}
\end{proposition}

\begin{proof}
The first part is clear and the second as $\H^n/\Pi$ is orientable
precisely when $\Pi$ acts orientation preservingly on $\H^n$, which is
to say, $\Pi$ is entirely composed of elements expressible as words
involving an even number of the generators (it is an elementary property
of Coxeter groups that the parity of any word is an invariant of the
corresponding element).\hfill\qed
\end{proof}

\section{Conjugacy in finite Coxeter groups}\label{conjugacy}



We saw in $\S 3$ that a list of representatives for the conjugacy classes of torsion 
elements in $(\Gamma,S)$ could be obtained by listing the conjugacy class representatives
of the finite subgroups generated by maximal $S'\subseteq S$. If
any such finite subgroup is reducible, 
then its conjugacy classes are easily obtained from those of its
irreducible components. Thus, we need only consider conjugacy in
irreducible finite Coxeter groups.

There are then two possibilities to consider: the Weyl groups and the
non-Weyl groups, and these are dealt with separately. 
Many of the Weyl groups have alternative descriptions in terms of
well-known finite groups. For example, those of type $A_n$ are
isomorphic to the
symmetric groups $\SS_{n+1}$; those of type $B_n$ arise as the split
extension
$$
1\longrightarrow\bigoplus^n\Z_2\longrightarrow B_n\longrightarrow
\SS_n\longrightarrow1,
$$
where for $\sigma\in\SS_n$ and $x_k$ the $\Z_2$-tuple with $1$ in the $k$-th
position and zeroes elsewhere, we have $\sigma^{-1}
x_k\sigma=x_{\sigma(k)}$ (the group of so-called signed permutations).
Similarly for those of type $D_n$ but with $\Z_2^n$ replaced by
$\Z_2^{n-1}$ (the even signed permutations).
Thus their conjugacy classes can be determined on a case by case basis.
Alternatively, \cite{Carter72} determines them conceptually using their
structure as
reflection groups, 
and it is the results from here that we use.

\subsection{The Weyl groups}\label{Carter}

The description is couched in terms of root systems 
(see \cite[Chapters 1 and 2]{Humphreys90}).
Let $V$ be a Euclidean space with orthogonal basis $\{\e_1,\ldots,\e_n\}$
and positive definite bilinear form defined by 
$\langle \e_i,\,\e_j\rangle=\delta_{ij}$, the Kronecker delta. For any
$\v\in V$, let $s_\v$ be the linear reflection in the hyperplane
orthogonal to $\v$. Tables \ref{roots1}-\ref{roots2} give the standard representations
of the Weyl groups acting on $V$: 
the subgroup of the orthogonal group $O(V)$ generated by the $s_\v$
where the $\v$ label the nodes of the symbol, is isomorphic to the
abstract Coxeter group having that symbol. These $\v$ are a set of simple roots.

If $(\Gamma,S)$ is an irreducible Weyl group, 
identify its symbol with the appropriate one in
Tables \ref{roots1}-\ref{roots2}, thus identifying $S$ with the $s_{\v}$.
It turns out \cite{Carter72} that the conjugacy classes in $\Gamma$ 
correspond to certain diagrams, closely related to 
the Coxeter symbol. We now proceed to give the
diagrams, and, with \S 3 in mind,
a method for obtaining a representative for the
conjugacy class in terms of the
generating reflections for $\Gamma$.


The diagrams for each group are given in Theorems \ref{Carter1}-\ref{Carter4} below. 
For each one, label the nodes (if they do not come prelabelled)
by roots from the appropriate root system in Tables \ref{roots1}--\ref{roots2}, 
in such a way that if the nodes labelled $\u$ and $\v$ are 
connected by an edge
labelled $m$, then
$$
\frac{4\langle \v,\u\rangle\langle \u,\v\rangle}
{\langle \v,\v\rangle\langle \u,\u\rangle}=m-2.
$$
In fact, we always have $m=2$ or $3$, and the diagrams use the same
labelling conventions as for Coxeter symbols. If the diagram can be
identified with a subsymbol of the Coxeter symbol, then a labelling is
easily obtained--take a labelling by simple roots off the Coxeter
symbol. 
Otherwise, it can be harder, as in the example below.

Colour the nodes of the labelled diagram black and white with a node
of one
colour joined only to nodes of the other. Let ${\mathscr B}$ and
${\mathscr W}$ be the
sets of nodes of the two colours. Then a representative of the conjugacy class
corresponding to the diagram is given by,
$$
\prod_{{\v}\in{\mathscr B}}s_{{\v}}\prod_{{\u}\in{\mathscr W}}s_{{\u}}.
$$
We use the algebraists convention of reading such expressions from left to right.
If $\v$ isn't a simple root, then $s_{\v}$ is not identified with one of the generating
reflections of $\Gamma$. In this case, find a simple $\v'$ and an
element $w\in\Gamma$, given in terms of the generators, such that 
the image of $\v$ under $w$ is $\v'$. Then $s_{\v}=ws_{\v'}w^{-1}$, an expression in
terms of the generating reflections.

\begin{table}
\begin{tabular}{lll}
\hline
Type and order&Root system&Coxeter symbol and simple system\\\hline
\begin{tabular}{l}
$A_n$\\
$(n+1)!$\\
\end{tabular}&$\{\e_i-\e_j \,\,(1\leq i\not= j\leq n+1)\}$&$\,\,\,\begin{picture}(150,30)
\put(-8,0){{$\e_1-\e_{2}$}}\put(25,17){{$\e_{2}-\e_{3}$}}
\put(75,0){{$\e_{n-1}-\e_{n}$}}\put(115,17){{$\e_{n}-\e_{n+1}$}}
\put(4,9){\circle{8}}\put(8,9){\line(1,0){30}}
\put(42,9){\circle{8}}\put(46,9){\line(1,0){12}}
\put(62,9){$\ldots$}\put(78,9){\line(1,0){12}}
\put(94,9){\circle{8}}\put(98,9){\line(1,0){30}}
\put(132,9){\circle{8}}                 
\end{picture}$\\
\begin{tabular}{l}
$D_n$\\
$2^{n-1}n!$\\
\end{tabular}&$\{\pm\e_i\pm\e_j\,\,(1\leq i<j\leq n)\}$&
\begin{picture}(150,55)
\put(-8,0){{$\e_1-\e_{2}$}}
\put(25,17){{$\e_{2}-\e_{3}$}}
\put(100,7){{$\e_{n-2}-\e_{n-1}$}}\put(110,42){{$\e_{n-1}-\e_{n}$}}
\put(110,-27){{$\e_{n-1}+\e_{n}$}}
\put(4,9){\circle{8}}\put(8,9){\line(1,0){30}}
\put(42,9){\circle{8}}\put(46,9){\line(1,0){12}}
\put(62,9){$\ldots$}\put(78,9){\line(1,0){12}}
\put(94,9){\circle{8}}                 
\put(96.83,11.83){\line(1,1){21.21}}    
\put(96.83,6.17){\line(1,-1){21.21}}   
\put(120.98,35.98){\circle{8}}          
\put(120.98,-17.98){\circle{8}}         
\end{picture}\\
\begin{tabular}{l}
$B_n$\\
$2^nn!$\\
\end{tabular}&
\begin{tabular}{l}
$\{\pm\e_i\,\,(1\leq i\leq n)$,\\
$\pm\e_i\pm\e_j\,\,(1\leq i<j\leq n)\}$\\
\end{tabular}
&\begin{picture}(150,50)
\put(-8,0){{$\e_1-\e_{2}$}}\put(25,17){{$\e_{2}-\e_3$}}
\put(75,0){{$\e_{n-1}-\e_{n}$}}\put(130,17){{$\e_{n}$}}
\put(109,11){$4$}
\put(4,9){\circle{8}}\put(8,9){\line(1,0){30}}
\put(42,9){\circle{8}}\put(46,9){\line(1,0){12}}
\put(62,9){$\ldots$}\put(78,9){\line(1,0){12}}
\put(94,9){\circle{8}}\put(97.5,9){\line(1,0){30}}
\put(132,9){\circle{8}}                 
\end{picture}\\\hline
\end{tabular}\caption{Root systems for the classical Weyl groups
\cite[\S 2.10]{Humphreys90}}\label{roots1}
\end{table}

For example, suppose we have a group of type $E_6$ with generators,
$$
\begin{picture}(150,40)
\put(2,10){{$x_1$}}\put(38,10){{$x_2$}}
\put(82,10){{$x_3$}}\put(116,10){{$x_4$}}
\put(86,36){{$x_6$}}\put(154,10){{$x_5$}}
\put(4,4){\circle{8}}\put(8,4){\line(1,0){30}}
\put(42,4){\circle{8}}\put(46,4){\line(1,0){30}}
\put(80,4){\circle{8}}\put(84,4){\line(1,0){30}}
\put(118,4){\circle{8}}\put(122,4){\line(1,0){30}}
\put(156,4){\circle{8}}
\put(80,36){\circle{8}}\put(80,8){\line(0,1){24}}
\end{picture}
$$
Consult the root system for $E_6$ in Table \ref{roots2} and make the 
identifications $x_1= s_{\e_5-\e_4}$, 
$x_2= s_{\e_4-\e_3}$, $x_3= s_{\e_3-\e_2}$, 
$x_4= s_{\e_2-\e_1}$, $x_5= s_{\v'}$ and 
$x_6= s_{\e_1+\e_2}$ where 
$\v'=\e_1+\e_8-\frac{1}{2}\sum_{i=1}^8 \e_i$.

By Theorem \ref{Carter3}, $E_6$ has a conjugacy class of elements of order three
corresponding to the diagram
$\coprod_{i=1}^3 \begin{picture}(32,6)
\put(3,3){\circle{6}}\put(6,3){\line(1,0){20}}
\put(29,3){\circle{6}}\end{picture}$. 
A labelling and a colouring is given by,
$$
%
\begin{picture}(180,20)
\put(26,18){{$\e_2-\e_1$}}
\put(-75,0){{$\e_1+\e_8-\frac{1}{2}(\sum_{i=1}^8 \e_i)$}}
\put(53,0){{$\e_4-\e_3$}}
\put(86,18){{$\e_5-\e_4$}}
\put(115,0){{$\e_1+\e_2$}}
\put(140,18){{$-\e_6-\e_7+\frac{1}{2}\sum_{i=1}^8 \e_i$}}





\put(8,10){\circle*{6}}\put(11,10){\line(1,0){30}}
\put(44,10){\circle{6}}
\put(68,10){\circle*{6}}\put(71,10){\line(1,0){30}}
\put(104,10){\circle{6}}
\put(128,10){\circle*{6}}\put(131,10){\line(1,0){30}}
\put(164,10){\circle{6}}\end{picture}
$$
There can be no labelling entirely by simple roots as the diagram is
not a subsymbol of the Coxeter symbol for $E_6$.
Thus, the roots labelling the nodes
of the coloured
diagram correspond, from left to right, to the reflections 
$x_5,x_4,x_2,x_1,x_6$
and $s_{\v}$ for $\v=-\e_6-\e_7+\frac{1}{2}\sum_{i=1}^8 \e_i$.
This last one is not a simple root, but if
$w=x_6x_3x_2x_1x_4x_3x_2x_6x_3x_4$ (reading from left to right), then
in the action of $\Gamma$ on $V$,
$\v$ is sent by $w$ to $\v'$ above.
When performing these calculations it is helpful to remember that the
reflection $s_{\e_i-\e_j}$ permutes the basis vectors by transposing
$\e_i$ and $\e_j$ while fixing all others.
Thus $s_{\v}=ws_{\v'}w^{-1}=wx_5w^{-1}$ and 
a representative of the conjugacy class corresponding to this
diagram is
$$
x_5x_2x_6x_4x_1wx_5w^{-1}.
$$

In \cite{Carter72} a product of this form is read from 
right to left, so the inverse of 
``our'' element is obtained (as the $s_{\v}$ are involutions). 
But all the diagrams below have the form
$$
\begin{picture}(88,6)
\put(-2,10){$\v_0$}
\put(24,10){$\v_1$}
\put(45,10){$\v_{k-1}$}
\put(76,10){$\v_k$}
\put(3,3){\circle{6}}\put(6,3){\line(1,0){20}}
\put(29,3){\circle{6}}\put(36,3){$\ldots$}
\put(55,3){\circle{6}}\put(58,3){\line(1,0){20}}
\put(81,3){\circle{6}}\end{picture}
,
$$ 
so that the subgroup generated by the $s_{\v_i}$ is 
isomorphic to $\SS_{k+1}$, and the two elements obtained are conjugate in $\Gamma$ anyway. 
Alternatively, although this is far less trivial, any element of a Weyl group is conjugate to its inverse
\cite[Corollary to Theorem C]{Carter72}.

\begin{table}
\begin{tabular}{ll}
\hline
Type, order and Root system&Coxeter symbol and simple system\\\hline
\begin{tabular}[t]{l}
$F_4$\\
$2^7\,3^2$.\\
$\{\pm\e_i\pm\e_j\,\,(1\leq i<j\leq 4),$\\
$\pm\e_i\,\,(1\leq i\leq 4),$\\
$\frac{1}{2}(\pm\e_1\pm\e_2\pm\e_3\pm\e_4)\};$\\
\end{tabular}
&$\,\,\,\begin{picture}(170,20)(20,20)
\put(-10,10){{$\e_1-\e_2$}}
\put(28,-6){{$\e_2-\e_3$}}
\put(76,10){{$\e_3$}}
\put(58,6){$4$}
\put(90,-6){{$\frac{1}{2}(-\e_1-\e_2-\e_3+\e_4)$}}
\put(4,4){\circle{8}}\put(8,4){\line(1,0){30}}
\put(42,4){\circle{8}}\put(46,4){\line(1,0){30}}
\put(80,4){\circle{8}}\put(84,4){\line(1,0){30}}
\put(118,4){\circle{8}}
\end{picture}$\\
\begin{tabular}{l}
$E_6$\\
$2^7\,3^4\,5$.\\
$\{\pm\e_i\pm\e_j\,\,(1\leq i<j\leq 5),$\\
$\frac{1}{2}\sum_{i=1}^8 \epsilon_i\e_i\};$\\
$(\epsilon_i=\pm 1, \prod\epsilon_i=1, \epsilon_8=-\epsilon_7=-\epsilon_6)$\\
\end{tabular}
&\begin{picture}(170,40)(35,30)
\put(-8,15){{$\e_5-\e_4$}}\put(30,30){{$\e_4-\e_3$}}
\put(68,15){{$\e_3-\e_2$}}\put(106,30){{$\e_2-\e_1$}}
\put(86,56){{$\e_1+\e_2$}}\put(130,15){{$\e_1+\e_8-\frac{1}{2}\sum_{i=1}^{8} \e_i$}}
\put(4,24){\circle{8}}\put(8,24){\line(1,0){30}}
\put(42,24){\circle{8}}\put(46,24){\line(1,0){30}}
\put(80,24){\circle{8}}\put(84,24){\line(1,0){30}}
\put(118,24){\circle{8}}\put(122,24){\line(1,0){30}}
\put(156,24){\circle{8}}
\put(80,56){\circle{8}}\put(80,28){\line(0,1){24}}
\end{picture}\\
\begin{tabular}{l}
$E_7$\\
$2^{10}\,3^4\,5\,7$.\\
$\{\pm\e_i\pm\e_j\,\,(1\leq i<j\leq 6),$\\
$\pm(\e_7-\e_8),$\\
$\frac{1}{2}\sum_{i=1}^8 \epsilon_i\e_i\};$\\
$(\epsilon_i=\pm 1, \prod\epsilon_i=1, \epsilon_8=-\epsilon_7)$\\\\
\end{tabular}&
\begin{picture}(170,60)(15,0)
\put(-32,-5){{$\e_6-\e_5$}}\put(6,-20){{$\e_5-\e_4$}}
\put(44,-5){{$\e_4-\e_3$}}\put(82,-20){{$\e_3-\e_2$}}
\put(120,-5){$\e_2-\e_1$}
\put(100,21){{$\e_1+\e_2$}}\put(144,-20){{$\e_1+\e_8-\frac{1}{2}\sum_{i=1}^{8} \e_i$}}
\put(-20,-11){\circle{8}}\put(-16,-11){\line(1,0){30}}
\put(18,-11){\circle{8}}\put(22,-11){\line(1,0){30}}
\put(56,-11){\circle{8}}\put(60,-11){\line(1,0){30}}
\put(94,-11){\circle{8}}\put(98,-11){\line(1,0){30}}
\put(132,-11){\circle{8}}\put(136,-11){\line(1,0){30}}
\put(170,-11){\circle{8}}
\put(94,21){\circle{8}}\put(94,-7){\line(0,1){24}}         
\end{picture}\\
\begin{tabular}{l}
$E_8$\\
$2^{14}\,3^5\,5^2\,7$\\\
$\{\pm\e_i\pm\e_j\,\,(1\leq i<j\leq 8),$\\
$\frac{1}{2}\sum_{i=1}^8 \epsilon_i\e_i\};$\\
$(\epsilon_i=\pm 1, \prod\epsilon_i=1)$\\\\
\end{tabular}&
\begin{picture}(170,40)(26,0)
\put(-42,-20){{$\e_7-\e_6$}}
\put(-2,-5){{$\e_6-\e_5$}}\put(36,-20){{$\e_5-\e_4$}}
\put(74,-5){{$\e_4-\e_3$}}\put(112,-20){{$\e_3-\e_2$}}
\put(150,-5){$\e_2-\e_1$}
\put(130,20){{$\e_1+\e_2$}}
\put(174,-20){{$\e_1+\e_8-\frac{1}{2}\sum_{i=1}^{8} \e_i$}}
\put(-30,-11){\circle{8}}\put(-26,-11){\line(1,0){30}}
\put(8,-11){\circle{8}}\put(12,-11){\line(1,0){30}}
\put(46,-11){\circle{8}}\put(50,-11){\line(1,0){30}}
\put(84,-11){\circle{8}}\put(88,-11){\line(1,0){30}}
\put(122,-11){\circle{8}}\put(126,-11){\line(1,0){30}}
\put(160,-11){\circle{8}}\put(164,-11){\line(1,0){30}}
\put(198,-11){\circle{8}}
\put(122,21){\circle{8}}\put(122,-7){\line(0,1){24}}         
\end{picture}\\\hline
\end{tabular}\caption{Root systems for the exceptional Weyl groups
\cite[\S 2.10]{Humphreys90}}\label{roots2}
\end{table}

\begin{theorem}[types $A$ and $B$]\label{Carter1}
The conjugacy classes of prime order in the Weyl groups of types $A_n$ and
$B_n$ correspond to the diagrams:
\begin{enumerate}
\item[$A_n$.] $\text{order }p\geq 2:$ $\displaystyle{\coprod^k 
\begin{picture}(88,6)
\put(0,0){$\underbrace{\vrule width 85pt height 0pt depth 0 pt}_{{p-1}}$}
\put(3,3){\circle{6}}\put(6,3){\line(1,0){20}}
\put(29,3){\circle{6}}\put(38,3){$\ldots$}
\put(55,3){\circle{6}}\put(58,3){\line(1,0){20}}
\put(81,3){\circle{6}}\end{picture}}$ for all $k\geq 1$ with $kp\leq n+1$.
\item[$B_n$.] 
$$\text{order }2: \coprod^k\begin{picture}(45,10)
\put(5,7){${\e}_i-{\e}_{i+1}$}
\put(20,1){\circle{6}}\end{picture}\coprod^m
\begin{picture}(20,10)
\put(8,7){${\e}_j$}
\put(10,1){\circle{6}}\end{picture}
\text{ and order }p\geq 3:
\coprod^k\begin{picture}(88,25)
\put(0,8){$\overbrace{\vrule width 85pt height 0pt depth 0 pt}^{{p-1}}$}
\put(3,3){\circle{6}}\put(6,3){\line(1,0){20}}
\put(29,3){\circle{6}}\put(38,3){$\ldots$}
\put(55,3){\circle{6}}\put(58,3){\line(1,0){20}}
\put(81,3){\circle{6}}\end{picture},
$$ 
where the order $2$ diagrams are for all $k+m> 0$ with $2k+m\leq n;$
the order $p\geq 3$ diagrams are for all $k\geq 1$ with $kp\leq n$.
\end{enumerate}
\end{theorem}

The theorem is Propositions 23 and 24 of \cite[\S 7]{Carter72}. Note that all the
unlabelled graphs above can be obtained as subsymbols of the Coxeter symbols for $A_n$ and
$B_n$.


\begin{theorem}[type $D$]\label{Carter2}
The conjugacy classes with prime order in $D_n$ correspond to the diagrams:
(a) order $2:$ $$\displaystyle{\coprod^k\begin{picture}(45,10)
\put(5,7){${\e}_i-{\e}_{i+1}$}
\put(20,1){\circle{6}}\end{picture}
\coprod^m\begin{picture}(32,10)
\put(5,7){${\e}_j-{\e}_{j+1}$}
\put(50,7){${\e}_j+{\e}_{j+1}$}
\put(20,1){\circle{6}}
\put(65,1){\circle{6}}\end{picture}}$$ 
for all $k+m> 0$ with $2(k+m)\leq n$,
except for when $n$ is even and $m=0$, where
$\coprod_{i=1}^{n/2} \begin{picture}(6,6)
\put(3,3){\circle{6}}\end{picture}$
corresponds to two classes with associated diagrams:
$$\begin{picture}(120,15)
\put(-12,10){$\e_1-\e_2$}\put(25,10){$\e_3-\e_4$}
\put(65,3){$\cdots$}
\put(83,10){$\e_{n-3}-\e_{n-2}$}\put(149,10){$\e_{n-1}\pm\e_n$}
\put(3,3){\circle{6}}
\put(40,3){\circle{6}}
\put(110,3){\circle{6}}
\put(167,3){\circle{6}}\end{picture}$$
(b) order $p\geq 3:$ $\displaystyle{\coprod^k\begin{picture}(88,25)
\put(0,8){$\overbrace{\vrule width 85pt height 0pt depth 0 pt}^{{p-1}}$}
\put(3,3){\circle{6}}\put(6,3){\line(1,0){20}}
\put(29,3){\circle{6}}\put(38,3){$\ldots$}
\put(55,3){\circle{6}}\put(58,3){\line(1,0){20}}
\put(81,3){\circle{6}}\end{picture}}$ for all $k\geq 1$ with $kp\leq n$.
\end{theorem}

See \cite[Proposition 25]{Carter72}. Note the diagram corresponding to two different conjugacy classes--hence the need for
different labellings. This exception arises from the fact that $D_n$ is isomorphic to
the subgroup of $B_n$ consisting of the even signed permutations of a set of size $n$: it permutes the basis 
$\{\e_1,\ldots,\e_n\}$ while changing the sign of an even number of them. From the proof
of Proposition 25 in \cite{Carter72}, if $w_1,w_2$ are the elements corresponding to the
two labellings in part (a) above, and $ww_1w^{-1}=w_2$, then $w$ must change the sign of an odd number of
the basis vectors, so the $w_i$ are not conjugate in $D_n$ although they
are in $B_n$. 

\begin{theorem}[type $E$]\label{Carter3}
The conjugacy classes with prime order in the Weyl groups of type 
$E_6,E_7$ and $E_8$ correspond to the diagrams:
\begin{enumerate}
\item[$E_6$.] 
(a) order $2:$ $\displaystyle{\coprod^k \begin{picture}(6,6)
\put(3,3){\circle{6}}\end{picture}, (k=1,\ldots, 4)};$ 
(b) order $3:$ $\displaystyle{\coprod^k \begin{picture}(32,15)
\put(3,3){\circle{6}}\put(6,3){\line(1,0){20}}
\put(29,3){\circle{6}}\end{picture}, (k=1,\ldots, 3)};$
(c) order $5:$ 
\begin{picture}(84,15)
\put(3,3){\circle{6}}\put(6,3){\line(1,0){20}}
\put(29,3){\circle{6}}\put(32,3){\line(1,0){20}}
\put(55,3){\circle{6}}\put(58,3){\line(1,0){20}}
\put(81,3){\circle{6}}\end{picture}.
\item[$E_7$.]
(a) order $2:$ $\displaystyle{\coprod^k \begin{picture}(6,6)
\put(3,3){\circle{6}}\end{picture}, (k=1,\ldots, 7)}$, except for $k=3,4$,
each of 
which correspond to two classes with associated diagrams:
$$
\begin{picture}(84,15)
\put(-12,10){$\e_1-\e_2$}\put(25,10){$\e_3-\e_4$}
\put(63,10){$\e_5\pm\e_6$}

\put(3,3){\circle{6}}
\put(40,3){\circle{6}}
\put(78,3){\circle{6}}\end{picture},
\,\,\,\,\,\,\,\,\,\,\,\,\,\,\,\,
\begin{picture}(130,15)
\put(-12,10){$\e_1-\e_2$}\put(25,10){$\e_3-\e_4$}
\put(63,10){$\e_5-\e_6$}\put(99,10){$\e_7-\e_8$}

\put(3,3){\circle{6}}
\put(40,3){\circle{6}}
\put(78,3){\circle{6}}
\put(114,3){\circle{6}}\end{picture}
$$
$$
\text{ and }\,\,\,\,\,\,\,\,
\begin{picture}(220,15)
\put(-12,10){$\e_1+\e_2$}\put(25,10){$\e_3+\e_4$}
\put(63,10){$\e_5+\e_6$}\put(122,1){$\e_1+\e_3+\e_5+\e_7-\frac{1}{2
}\sum_{i=1}^8\e_i$}

\put(3,3){\circle{6}}
\put(40,3){\circle{6}}
\put(78,3){\circle{6}}
\put(114,3){\circle{6}}\end{picture}
$$
(b) order $3:$ $\displaystyle{\coprod^k \begin{picture}(35,15)
\put(3,3){\circle{6}}\put(6,3){\line(1,0){20}}
\put(29,3){\circle{6}}\end{picture}, (k=1,\ldots, 3)};$
(c) order $5:$ 
\begin{picture}(86,15)
\put(3,3){\circle{6}}\put(6,3){\line(1,0){20}}
\put(29,3){\circle{6}}\put(32,3){\line(1,0){20}}
\put(55,3){\circle{6}}\put(58,3){\line(1,0){20}}
\put(81,3){\circle{6}}\end{picture},
(d) order $7:$ 
\begin{picture}(138,15)
\put(3,3){\circle{6}}\put(6,3){\line(1,0){20}}
\put(29,3){\circle{6}}\put(32,3){\line(1,0){20}}
\put(55,3){\circle{6}}\put(58,3){\line(1,0){20}}
\put(81,3){\circle{6}}\put(84,3){\line(1,0){20}}
\put(107,3){\circle{6}}\put(110,3){\line(1,0){20}}
\put(133,3){\circle{6}}\end{picture}.
\item[$E_8$.]
(a) order $2:$ $\displaystyle{\coprod^k \begin{picture}(6,6)
\put(3,3){\circle{6}}\end{picture}, (k=1,\ldots, 8)}$, except for $k=4$
which corresponds to two classes with associated diagrams:
$$
\begin{picture}(130,15)
\put(-12,10){$\e_1-\e_2$}\put(25,10){$\e_3-\e_4$}
\put(63,10){$\e_5-\e_6$}\put(105,10){$\e_7\pm\e_8$}

\put(3,3){\circle{6}}
\put(40,3){\circle{6}}
\put(78,3){\circle{6}}
\put(120,3){\circle{6}}\end{picture};
$$
(b) order $3:$ $\displaystyle{\coprod^k \begin{picture}(32,15)
\put(3,3){\circle{6}}\put(6,3){\line(1,0){20}}
\put(29,3){\circle{6}}\end{picture}, (k=1,\ldots, 4)};$
(c) order $5:$ $\displaystyle{\coprod^k
\begin{picture}(84,15)
\put(3,3){\circle{6}}\put(6,3){\line(1,0){20}}
\put(29,3){\circle{6}}\put(32,3){\line(1,0){20}}
\put(55,3){\circle{6}}\put(58,3){\line(1,0){20}}
\put(81,3){\circle{6}}\end{picture}, (k=1,2)};$
(d) order $7:$
\begin{picture}(140,15)
\put(3,3){\circle{6}}\put(6,3){\line(1,0){20}}
\put(29,3){\circle{6}}\put(32,3){\line(1,0){20}}
\put(55,3){\circle{6}}\put(58,3){\line(1,0){20}}
\put(81,3){\circle{6}}\put(84,3){\line(1,0){20}}
\put(107,3){\circle{6}}\put(110,3){\line(1,0){20}}
\put(133,3){\circle{6}}\end{picture}.
\end{enumerate}
\end{theorem}

Both Theorem \ref{Carter3} and Theorem \ref{Carter4} below arise from a careful examination of the 
results from \cite[\S 8]{Carter72}. These groups do not have convenient descriptions
as the classical groups do. In \cite{Carter72}, a collection of non-conjugate elements is
formed, and a counting argument gives that the list is complete. As with the groups of type
$D_n$, there are some pairs of conjugacy classes corresponding to the same (unlabelled)
diagram. That the
different labellings given in the two theorems yield non-conjugate elements is checked by
observing that they fix a different number of roots in the systems from
Table \ref{roots2}.

\begin{theorem}[type $F$]\label{Carter4}
The conjugacy classes with prime order in the Weyl group of type $F_4$ correspond to the diagrams:
(a) order $2:$ 
$$\begin{picture}(10,15)
\put(-12,10){$\e_1-\e_2$}
\put(3,3){\circle{6}}\end{picture},\,\,\,\,\,\,\,\,
\begin{picture}(10,15)
\put(-1,10){$\e_1$}
\put(3,3){\circle{6}}\end{picture},\,\,\,\,\,\,\,\,
\begin{picture}(45,15)
\put(-12,10){$\e_1-\e_2$}\put(25,10){$\e_3-\e_4$}
\put(3,3){\circle{6}}
\put(40,3){\circle{6}}\end{picture},
$$
$$
\begin{picture}(45,15)
\put(-12,10){$\e_1-\e_2$}\put(36,10){$\e_3$}
\put(3,3){\circle{6}}
\put(40,3){\circle{6}}\end{picture},\,\,\,\,\,\,\,\,
\begin{picture}(84,15)
\put(-12,-8){$\e_1-\e_2$}\put(25,-8){$\e_3-\e_4$}
\put(63,-8){$\e_3+\e_4$}
\put(3,3){\circle{6}}
\put(40,3){\circle{6}}
\put(78,3){\circle{6}}\end{picture}
\,\,\,\,\,\,\,\,\text{ and }\,\,\,\,\,\,\,\,
\begin{picture}(84,15)
\put(-12,10){$\e_1-\e_2$}\put(25,10){$\e_3-\e_4$}
\put(63,10){$\frac{1}{2}\sum_{i=1}^4\e_i$}

\put(3,3){\circle{6}}
\put(40,3){\circle{6}}
\put(78,3){\circle{6}}\end{picture};
$$
(b) order $3:$ 
$$
\begin{picture}(42,15)
\put(-14,-8){$\e_1-\e_2$}\put(23,-8){$\e_2-\e_3$}

\put(3,3){\circle{6}}\put(6,3){\line(1,0){30}}
\put(39,3){\circle{6}}\end{picture},\,\,\,\,\,\,\,\,
\begin{picture}(42,15)
\put(-1,10){$\e_4$}\put(23,10){$\frac{1}{2}\sum_{i=1}^4\e_i$}

\put(3,3){\circle{6}}\put(6,3){\line(1,0){30}}
\put(39,3){\circle{6}}\end{picture},\,\,\,\,\,\,\,\,\,\,\,\,\,\,\,\,\,\,\,\,\,\,
\begin{picture}(100,15)
\put(-14,-8){$\e_1-\e_2$}\put(23,-8){$\e_2-\e_3$}
\put(59,10){$\e_4$}\put(84,10){$\frac{1}{2}\sum_{i=1}^4\e_i$}

\put(3,3){\circle{6}}\put(6,3){\line(1,0){30}}
\put(39,3){\circle{6}}\put(63,3){\circle{6}}\put(66,3){\line(1,0){30}}
\put(99,3){\circle{6}}\end{picture}
$$
\end{theorem}

\begin{center}
\begin{table}
\begin{tabular}{lll}
\hline
Type and order&&Coxeter symbol and simple system\\\hline
\begin{tabular}[t]{l}
$H_3$\\
$2^3\,3\,5$.\\
\end{tabular}&
&$\,\,\,\begin{picture}(170,20)(0,5)
\put(23,10){{$\frac{1}{2}+b\i-a\j$}}
\put(59,-8){{$-a+\frac{1}{2}\i+b\j$}}
\put(96,6){$5$}
\put(124,2){{$a-\frac{1}{2}\i+b\j$}}
\put(42,4){\circle{8}}\put(46,4){\line(1,0){30}}
\put(80,4){\circle{8}}\put(84,4){\line(1,0){30}}
\put(118,4){\circle{8}}
\end{picture}$\\
\begin{tabular}[t]{l}
$H_4$\\
$2^6\,3^2\,5^2$.\\
\end{tabular}&
&$\,\,\,\begin{picture}(170,20)(0,5)
\put(-20,-8){{$-\frac{1}{2}-a\i+b\k$}}
\put(23,10){{$\frac{1}{2}+b\i-a\j$}}
\put(59,-8){{$-a+\frac{1}{2}\i+b\j$}}
\put(96,6){$5$}
\put(124,2){{$a-\frac{1}{2}\i+b\j$}}
\put(4,4){\circle{8}}\put(8,4){\line(1,0){30}}
\put(42,4){\circle{8}}\put(46,4){\line(1,0){30}}
\put(80,4){\circle{8}}\put(84,4){\line(1,0){30}}
\put(118,4){\circle{8}}
\end{picture}$\\
\\
\hline
\end{tabular}\caption{Root systems for type $H$ non-Weyl groups
\cite[\S 2.13]{Humphreys90}}\label{roots3}
\end{table}
\end{center}

\subsection{The non-Weyl groups}\label{}

Identify a $4$-dimensional Euclidean space $V$ with the division ring
$\HH$ of quaternions, so that the bilinear form of \S 4.1 becomes
$\langle \u,\v\rangle=\frac{1}{2}(\u\bar{\v}+\v\bar{\u})$ where
$\bar{\v}=c_1-c_2\i-c_3\j-c_4\k$ is the usual quaternionic conjugation.
If $\v\in\HH$ has norm $1$ then the reflection $s_{\v}$ is given by
$s_{\v}(\u)=-\u\bar{\v}\u$. Let the split extension of $\Z_2^4$ by
$\AA_4$ act on $\HH$ with the alternating group $\AA_4$ permuting
coordinates and the $\Z_2^4$ generated by sign changes. If 
$$
a=\frac{1+\sqrt{5}}{4}\text{ and }b=\frac{-1+\sqrt{5}}{4},
$$
then the $120$ images under this action of $1,\frac{1}{2}(1+\i+\j+\k)$
and $a+\frac{1}{2}\i+b\j$ form a root system for $H_4$ (\cite[\S
2.13]{Humphreys90}). The $30$ roots orthogonal to $\k$ give a root
system for $H_3$. Coxeter symbols and simple systems are given in Table
\ref{roots3}. 

\begin{proposition}[type $H$]\label{nonweyl1}
The conjugacy classes with prime order in the groups of type $H_3$ and
$H_4$ have representatives:
\begin{enumerate}
\item[$H_3$.] 
(a) order $2:$ $x_1,x_1x_3$ and $(x_2x_1x_3)^5;$ (b) order $3:$ $x_1x_2;$
(c) order $5:$ $x_2x_3$ and $(x_2x_3)^2$.
\item[$H_4$.] 
(a) order $2:$ $x_1, x_1x_4, x_2x_1x_3x_2x_1w^{12}\bar{w}x_4x_3x_4$ and
$x_1x_2x_1x_3x_2x_1w^{12}\bar{w}x_4x_3x_4;$ 
(b) order $3:$ $x_1x_2$  and $x_2x_3;$
(c) order $5:$ $x_3x_4$, $(x_3x_4)^2$, 
$x_3w^3\bar{w}w^2$, $x_3w^9\bar{w}w^2$ and
$x_1x_2x_1x_3w^8\bar{w}x_4$ $x_3x_4$,
\end{enumerate}
where $x_i=s_{\v_i}$ with $\v_i$ the $i$-th simple root from the left in
the Coxeter symbol, $\bar{w}=x_4x_3x_2$ and $w$ is the Coxeter element $x_4x_3x_2x_1$.
\end{proposition}

\begin{proof}
The $34$ conjugacy classes in $H_4$, and a representative transformation
of $\HH$ in each, are given in \cite[Table 3]{Grove74}. By observing
their effect on $1,\i,\j,\k$ the order of these transformations can be
computed to give four of order $2$, two of order $3$ and five of order
$5$, the others non-prime. Thus it remains to show that the elements stated in the Proposition
are non-conjugate.
For the order $2$ and $3$ elements this is most easily done by computing
the number of roots fixed by each--a tedious but finite task that gives for
example in the order $3$ case, two roots fixed by $x_1x_2$ and six by $x_2x_3$, Three of the
order five elements fix no roots, but nevertheless have different
traces. 

Similar calculations give the $H_3$ classes.
Alternatively, the group is well known to be isomorphic to $\Z_2\times\AA_5$,
with 
$\AA_5$ the rotations of a dodecahedron, and the center $\Z_2$
generated by the antipodal map. 
\hfill\qed
\end{proof}

Finally we have the, obviously well known,

\begin{proposition}[type $I$]\label{nonweyl2}
The conjugacy classes with prime order in the dihedral groups, 
$$
\,\,\,\begin{picture}(60,10)
\put(19,6){$m$}\put(0,-6){$x_1$}\put(38,-6){$x_2$}
\put(4,4){\circle{8}}\put(8,4){\line(1,0){30}}
\put(42,4){\circle{8}}
\end{picture}
$$
of type $I$ have representatives:
for $m=2k+1$, $x_1$ and $(x_1x_2)^l$; and for $m=2k$, $x_1,x_2$ and $(x_1x_2)^l$;
for all $l\leq k$ with $m/\gcd(m,l)$ is prime.
\end{proposition}





\section{Volume}\label{volume}






Let $\Pi$ be a group acting properly, freely, cofinitely by isometries on $\H^n$
and let $M=\H^n/\Pi$. If $n$ is even then the Gauss-Bonnet (-Hirzebruch)
formula \cite{Spivak79,Gromov82,Hirzebruch57,Hirzebruch66}, gives
$$
\vol(M)=\kappa_n \chi_{\top}(M),
$$
with $\chi_{\top}(M)$ the Euler characteristic of $M$ and 
$\kappa_n=2^n (n!)^{-1} (-\pi)^{n/2} (n/2)!$ As $\Pi$ is of finite
homological type and torsion free
the Euler characteristic, $\chi(\Pi)=\sum_i
(-1)^i\text{rank}_{\Z}H_i(\Pi)$ is defined \cite[\S IX.6]{Brown82}, and 
since $M$ is a
$K(\Pi,1)$ space ($\H^n$ being contractible) the homologies $H_*(M)\cong
H_*(\Pi)$, so $\chi_{\top}(M)=\chi(\Pi)$.

Now suppose $\Pi$ arises via a torsion free transitive $\C\Gamma$-module
$U$ as in \S 3, with $\Gamma$ a hyperbolic Coxeter group and $\dim U=m$
(which is thus the index in $\Gamma$ of $\Pi$).
Then $\chi(\Gamma)$ is defined and $\chi(\Pi)=m\chi(\Gamma)$. But calculating the Euler
characteristic of Coxeter groups turns out to be a simple task:

\begin{theorem}\label{Chiswell}
Let $(\Gamma,S)$ be a Coxeter group. 
For any $\ss\in\FF$, let $|\ss|$ be the order of this group. Then,
\begin{equation}\label{eulerchar}
\chi(\Gamma)=
\sum_{\substack{\ss_0<\cdots<\ss_k\\\ss_i\in\FF}}\frac{(-1)^k}{|\ss_0|},
\end{equation}
where the sum is over all chains $\ss_0<\cdots<\ss_k$ with $\ss_i\in\FF$.
\end{theorem}

A similar result appears in \cite{Chiswell92} for the Euler characteristic defined in \cite{Chiswell76}.

\begin{proof}
Let $\Sigma$ be the affine $\Gamma$-complex from \S 2, metrised as
there to be a CAT(0) space. For each $k$-cell $\ss=(\ss_0<\ss_1<\cdots
<\ss_k)$ let $\Gamma_{\ss}$ be the isotropy group of $\ss$ and $Y$ the
$k$-cells $\ss_0<\ss_1<\cdots <\ss_k$ where each $\ss_i\in\FF$. Then 
$Y$ is a set  of representatives of the cells of
$\Sigma$ modulo the $\Gamma$-action. Moreover, $\Sigma$ is contractible
(as indeed is any CAT(0) space) so by \cite[Proposition 7.3(e$^{\prime}$)]{Brown82},
$$
\chi(\Gamma)=\chi_{\Gamma}(\Sigma)=\sum_{\ss\in Y} (-1)^{\dim\ss} \chi(\Gamma_{\ss}),
$$
with the right hand side the $\Gamma$-equivariant Euler characteristic of $\Sigma$. 
If $\ss=(\ss_0<\cdots<\ss_k)\in Y$ then $\Gamma_\ss=\ss_0$ and $\chi(\Gamma_\ss)=1/|\Gamma_\ss|=1/|\ss_0|$ 
as these groups are finite.\hfill\qed
\end{proof}




If $\Gamma$ is a Coxeter symbol and $\Delta$ a subsymbol, let
$\Gamma\setminus\Delta$ be the subsymbol obtained by removing from
$\Gamma$ the vertices of $\Delta$ and their incident edges. 
The following lemma is then sometimes useful.

\begin{proposition}\label{eulerchar2}
Let $\Gamma$ be a Coxeter symbol and $\Psi$ a subsymbol.
Let $\Sigma_{\Psi}$ be the sum (\ref{eulerchar})
restricted to those chains $\ss_0<\cdots<\ss_k$ with $\Psi<\ss_0$.
Then,
$$
\Sigma_{\Psi}=\sum_{\Delta\in{\mathscr P}(\Psi)} (-1)^{v(\Delta)}
\chi(\Gamma\setminus\Delta),
$$
where $\PP(\Psi)$ is the set of all subsymbols of $\Psi$, and $v(\Delta)$ the number of nodes
of $\Delta$.
\end{proposition}

\begin{proof}
Apply the inclusion-exclusion principle as in \cite[\S 2.1]{Stanley97}\hfill\qed
\end{proof}


If $\Psi$ is the symbol for an infinite Coxeter group, then $\sum_{\Psi}=0$, a sum
over an empty set, so in particular if $\Gamma$ is infinite we obtain as a Corollary that,
$$
\sum_{\Delta\in\PP(\Gamma)}(-1)^{v(\Delta)}\chi(\Gamma\setminus\Delta)=0,
$$
which is Serre's closed form for the Euler characteristic \cite[Proposition 16]{Serre71}

For $n$ odd, the Euler characteristic of an $n$-dimensional hyperbolic Coxeter group
is zero,
and there is in general no simple method for computing volumes. Of
course, $\vol(\H^n/\Pi)=m\vol(P)$ where $\Pi$ is an index $m$ subgroup
of Coxeter group $\Gamma$ with fundamental polytope $P$. For certain
$\Gamma$ the volumes of such $P$ has been determined \cite{Ratcliffe99}.

\section{Examples}\label{examples}

We give three manifolds as examples of the construction, one in each $\H^n$ for $n=4,5,6$. 
We have let volume be our guiding principle, so that the three are either the smallest possible 
($n=4$) or would appear to be the smallest known ($n=5,6$). As mentioned in the Introduction, a large
family of $4$-manifolds with $\chi=1$ have recently been constructed by Ratcliffe and Tschantz 
\cite{Ratcliffe00} using a computer.

We first explain the notation. In each case we have a Coxeter group $\Gamma$ acting on a collection 
of finite sets $\Omega_i$, so
that if $U_i$ is the $\C$-vector space with basis $\Omega_i$, then
$\otimes U_i$ is a torsion free $\C\Gamma$-module.

The actions of $\Gamma$ on the $\om_i$ are depicted using two notations.
If $\om_i$ has solid black nodes (as in 6.1, 6.2 diagrams (1)--(3) and
6.3 diagrams (1)--(2)) then generator $x_j$ of $\Gamma$ acts by swapping two nodes
if they are connected by an edge labelled 
according to the scheme:
$\setlength{\unitlength}{0.00083333in}
\begingroup\makeatletter\ifx\SetFigFont\undefined%
\gdef\SetFigFont#1#2#3#4#5{%
  \reset@font\fontsize{#1}{#2pt}%
  \fontfamily{#3}\fontseries{#4}\fontshape{#5}%
  \selectfont}%
\fi\endgroup%
{\renewcommand{\dashlinestretch}{30}
\begin{picture}(450,50)(0,0)
\path(0,50)(450,50)
\end{picture}}$ 
for $x_1$;
$\setlength{\unitlength}{0.00083333in}
\begingroup\makeatletter\ifx\SetFigFont\undefined%
\gdef\SetFigFont#1#2#3#4#5{%
  \reset@font\fontsize{#1}{#2pt}%
  \fontfamily{#3}\fontseries{#4}\fontshape{#5}%
  \selectfont}%
\fi\endgroup%
{\renewcommand{\dashlinestretch}{30}
\begin{picture}(450,50)(0,0)
\path(0,40)(450,40)
\path(0,60)(450,60)
\end{picture}}$
for $x_2$;
$\setlength{\unitlength}{0.00083333in}
\begingroup\makeatletter\ifx\SetFigFont\undefined%
\gdef\SetFigFont#1#2#3#4#5{%
  \reset@font\fontsize{#1}{#2pt}%
  \fontfamily{#3}\fontseries{#4}\fontshape{#5}%
  \selectfont}%
\fi\endgroup%
{\renewcommand{\dashlinestretch}{30}
\begin{picture}(450,50)(0,0)
\path(0,50)(450,50)
\put(225,50){\ellipse{70}{70}}
\put(225,50){\whiten\ellipse{70}{70}}
\put(225,50){\blacken\ellipse{30}{30}}
\put(225,50){\ellipse{30}{30}}
\end{picture}}$
for $x_3$;
$\setlength{\unitlength}{0.00083333in}
\begingroup\makeatletter\ifx\SetFigFont\undefined%
\gdef\SetFigFont#1#2#3#4#5{%
  \reset@font\fontsize{#1}{#2pt}%
  \fontfamily{#3}\fontseries{#4}\fontshape{#5}%
  \selectfont}%
\fi\endgroup%
{\renewcommand{\dashlinestretch}{30}
\begin{picture}(450,50)(0,0)
\path(0,50)(450,50)
\put(150,50){\whiten\ellipse{50}{50}}
\put(150,50){\ellipse{50}{50}}
\put(300,50){\whiten\ellipse{50}{50}}
\put(300,50){\ellipse{50}{50}}
\end{picture}}$
for $x_4$;
$\setlength{\unitlength}{0.00083333in}
\begingroup\makeatletter\ifx\SetFigFont\undefined%
\gdef\SetFigFont#1#2#3#4#5{%
  \reset@font\fontsize{#1}{#2pt}%
  \fontfamily{#3}\fontseries{#4}\fontshape{#5}%
  \selectfont}%
\fi\endgroup%
{\renewcommand{\dashlinestretch}{30}
\begin{picture}(450,50)(0,0)
\path(0,50)(450,50)
\put(150,50){\blacken\ellipse{40}{40}}
\put(150,50){\ellipse{40}{40}}
\put(300,50){\blacken\ellipse{40}{40}}
\put(300,50){\ellipse{40}{40}}
\end{picture}}$
for $x_5$;
$\setlength{\unitlength}{0.00083333in}
\begingroup\makeatletter\ifx\SetFigFont\undefined%
\gdef\SetFigFont#1#2#3#4#5{%
  \reset@font\fontsize{#1}{#2pt}%
  \fontfamily{#3}\fontseries{#4}\fontshape{#5}%
  \selectfont}%
\fi\endgroup%
{\renewcommand{\dashlinestretch}{30}
\begin{picture}(450,50)(0,0)
\path(0,50)(450,50)
\put(225,50){\blacken\ellipse{40}{40}}
\put(225,50){\ellipse{40}{40}}
\end{picture}}$
for $x_6$
and 
$\setlength{\unitlength}{0.00083333in}
\begingroup\makeatletter\ifx\SetFigFont\undefined%
\gdef\SetFigFont#1#2#3#4#5{%
  \reset@font\fontsize{#1}{#2pt}%
  \fontfamily{#3}\fontseries{#4}\fontshape{#5}%
  \selectfont}%
\fi\endgroup%
{\renewcommand{\dashlinestretch}{30}
\begin{picture}(450,50)(0,0)
\path(0,50)(450,50)
\put(225,50){\whiten\ellipse{50}{50}}
\put(225,50){\ellipse{50}{50}}
\end{picture}}$ for $x_7$ (or even combinations that can be interpreted unambiguously like
$\setlength{\unitlength}{0.00083333in}
\begingroup\makeatletter\ifx\SetFigFont\undefined%
\gdef\SetFigFont#1#2#3#4#5{%
  \reset@font\fontsize{#1}{#2pt}%
  \fontfamily{#3}\fontseries{#4}\fontshape{#5}%
  \selectfont}%
\fi\endgroup%
{\renewcommand{\dashlinestretch}{30}
\begin{picture}(450,50)(0,0)
\path(0,40)(450,40)
\path(0,60)(450,60)
\put(225,50){\ellipse{70}{70}}
\put(225,50){\whiten\ellipse{70}{70}}
\put(225,50){\blacken\ellipse{30}{30}}
\put(225,50){\ellipse{30}{30}}
\end{picture}}$ for an $x_2$ and $x_3$ edge). Nodes not so connected are fixed by the 
corresponding $x_j$.

The other notation, used in \S 6.2 diagram $\om_4$ and \S 6.3 diagrams
$\om_3$ and $\om_4$, is convenient when
the action of $\Gamma$ is imprimitive and we want to compress rather a
lot of information (so the resulting diagrams
require a certain amount of
decoding). 
The numbered vertices in
the diagrams are the blocks of imprimitivity of an imprimitive
$\Gamma$-action; the action on the blocks is depicted 
using the notation above.

To recover the action on the original set $\om$, suppose that the block 
$\om_i=\{m_{i1},m_{i2},\ldots,m_{ik}\}$. 
We
use the pair
%
$(i,j)_k=\ss$
to mean $m_{il}$ goes to $m_{j\sigma(l)}$ under the action of the
generator $x_k$.
An absence of such an indication corresponds to $\sigma$ being the identity.
If the block $\om_i$ is fixed by $x_j$, then the action of $x_j$ on the $m_{il}$ is
described with
each diagram.

\subsection{A $4$-manifold with $\chi=1$.}

Let $\Gamma$ be the Coxeter group with symbol on the left,
$$
\begin{tabular}{ccc}
\begin{picture}(84,88)(0,-5)
\put(88,40){${x_2}$}\put(50,78){${x_1}$}\put(50,2){${x_3}$}
\put(-14,40){${x_4}$}\put(46,46){${x_5}$}
\put(20,44){$4$}
\put(4,42){\circle{8}}\put(8,42){\line(1,0){30}}
\put(42,42){\circle{8}}\put(46,42){\line(1,0){30}}
\put(80,42){\circle{8}}
\put(42,46){\line(0,1){30}}\put(42,38){\line(0,-1){30}}
\put(42,80){\circle{8}}\put(42,4){\circle{8}}\end{picture}
&
\begin{tabular}[b]{c}
\begin{tabular}[b]{c}
\setlength{\unitlength}{0.00083333in}
\begingroup\makeatletter\ifx\SetFigFont\undefined%
\gdef\SetFigFont#1#2#3#4#5{%
  \reset@font\fontsize{#1}{#2pt}%
  \fontfamily{#3}\fontseries{#4}\fontshape{#5}%
  \selectfont}%
\fi\endgroup%
{\renewcommand{\dashlinestretch}{30}
\begin{picture}(1016,333)(0,-10)
\put(733.000,65.250){\arc{487.500}{3.5364}{5.8884}}
\put(733.000,252.750){\arc{487.500}{0.3948}{2.7468}}
\put(508,159){\blacken\ellipse{100}{100}}
\put(508,159){\ellipse{100}{100}}
\put(958,159){\blacken\ellipse{100}{100}}
\put(958,159){\ellipse{100}{100}}
\put(58,159){\blacken\ellipse{100}{100}}
\put(58,159){\ellipse{100}{100}}

\path(508,169)(958,169)
\path(508,149)(958,149)
	
\path(58,159)(508,159)

\put(733,15){\ellipse{70}{70}}
\put(733,15){\whiten\ellipse{70}{70}}
\put(733,15){\blacken\ellipse{30}{30}}
\put(733,15){\ellipse{30}{30}}

\put(200,159){\blacken\ellipse{40}{40}}
\put(200,159){\ellipse{40}{40}}
\put(350,159){\blacken\ellipse{40}{40}}
\put(350,159){\ellipse{40}{40}}

\end{picture}
}\\
$\Omega_1$\\
\end{tabular}
\\
\\
%







\begin{tabular}[b]{c}
\setlength{\unitlength}{0.00083333in}
\begingroup\makeatletter\ifx\SetFigFont\undefined%
\gdef\SetFigFont#1#2#3#4#5{%
  \reset@font\fontsize{#1}{#2pt}%
  \fontfamily{#3}\fontseries{#4}\fontshape{#5}%
  \selectfont}%
\fi\endgroup%
{\renewcommand{\dashlinestretch}{30}
\begin{picture}(1668,579)(0,-10)
\put(65.250,282.000){\arc{487.500}{5.1072}{7.4592}}
\put(515.250,282.000){\arc{487.500}{5.1072}{7.4592}}
\put(1415.250,282.000){\arc{487.500}{5.1072}{7.4592}}
\put(252.750,282.000){\arc{487.500}{1.9656}{4.3176}}

\put(1059,507){\blacken\ellipse{100}{100}}
\put(1059,507){\ellipse{100}{100}}
\put(1509,507){\blacken\ellipse{100}{100}}
\put(1509,507){\ellipse{100}{100}}
\put(609,507){\blacken\ellipse{100}{100}}
\put(609,507){\ellipse{100}{100}}
\put(159,507){\blacken\ellipse{100}{100}}
\put(159,507){\ellipse{100}{100}}
\put(159,57){\blacken\ellipse{100}{100}}
\put(159,57){\ellipse{100}{100}}
\put(1059,57){\blacken\ellipse{100}{100}}
\put(1059,57){\ellipse{100}{100}}
\put(1509,57){\blacken\ellipse{100}{100}}
\put(1509,57){\ellipse{100}{100}}
\put(609,57){\blacken\ellipse{100}{100}}
\put(609,57){\ellipse{100}{100}}

\path(1059,507)(1509,507)
\path(609,507)(1059,507)
\path(159,57)(159,507)
\path(609,57)(609,507)

\path(1069,57)(1069,507)
\path(1049,57)(1049,507)

\path(1519,57)(1519,507)
\path(1499,57)(1499,507)

\path(609,57)(1059,57)

\path(159,67)(609,67)
\path(159,47)(609,47)

\path(1059,57)(1509,57)

\path(159,517)(609,517)
\path(159,497)(609,497)

\path(1059,507)(1509,57)
\path(1509,507)(1059,57)

\put(159,282){\ellipse{70}{70}}
\put(159,282){\whiten\ellipse{70}{70}}
\put(159,282){\blacken\ellipse{30}{30}}
\put(159,282){\ellipse{30}{30}}

\put(609,282){\ellipse{70}{70}}
\put(609,282){\whiten\ellipse{70}{70}}
\put(609,282){\blacken\ellipse{30}{30}}
\put(609,282){\ellipse{30}{30}}

\put(1284,507){\ellipse{70}{70}}
\put(1284,507){\whiten\ellipse{70}{70}}
\put(1284,507){\blacken\ellipse{30}{30}}
\put(1284,507){\ellipse{30}{30}}

\put(1284,57){\ellipse{70}{70}}
\put(1284,57){\whiten\ellipse{70}{70}}
\put(1284,57){\blacken\ellipse{30}{30}}
\put(1284,57){\ellipse{30}{30}}

\put(759,507){\blacken\ellipse{40}{40}}
\put(759,507){\ellipse{40}{40}}
\put(909,507){\blacken\ellipse{40}{40}}
\put(909,507){\ellipse{40}{40}}

\put(759,57){\blacken\ellipse{40}{40}}
\put(759,57){\ellipse{40}{40}}
\put(909,57){\blacken\ellipse{40}{40}}
\put(909,57){\ellipse{40}{40}}

\put(159,175){\whiten\ellipse{50}{50}}
\put(159,175){\ellipse{50}{50}}
\put(159,389){\whiten\ellipse{50}{50}}
\put(159,389){\ellipse{50}{50}}

\put(609,175){\whiten\ellipse{50}{50}}
\put(609,175){\ellipse{50}{50}}
\put(609,389){\whiten\ellipse{50}{50}}
\put(609,389){\ellipse{50}{50}}

\put(1059,207){\whiten\ellipse{50}{50}}
\put(1059,207){\ellipse{50}{50}}
\put(1059,357){\whiten\ellipse{50}{50}}
\put(1059,357){\ellipse{50}{50}}

\put(1509,207){\whiten\ellipse{50}{50}}
\put(1509,207){\ellipse{50}{50}}
\put(1509,357){\whiten\ellipse{50}{50}}
\put(1509,357){\ellipse{50}{50}}

\put(20,207){\blacken\ellipse{40}{40}}
\put(20,207){\ellipse{40}{40}}
\put(20,357){\blacken\ellipse{40}{40}}
\put(20,357){\ellipse{40}{40}}

\put(1648,207){\blacken\ellipse{40}{40}}
\put(1648,207){\ellipse{40}{40}}
\put(1648,357){\blacken\ellipse{40}{40}}
\put(1648,357){\ellipse{40}{40}}

\end{picture}
}\\
$\Omega_2$\\
\end{tabular}
\end{tabular}
&
\begin{tabular}[b]{c}
\setlength{\unitlength}{0.00083333in}
\begingroup\makeatletter\ifx\SetFigFont\undefined%
\gdef\SetFigFont#1#2#3#4#5{%
  \reset@font\fontsize{#1}{#2pt}%
  \fontfamily{#3}\fontseries{#4}\fontshape{#5}%
  \selectfont}%
\fi\endgroup%
{\renewcommand{\dashlinestretch}{30}
\begin{picture}(1220,1737)(0,-10)
\put(610.000,759.000){\arc{1500.000}{0.9273}{2.2143}}
\put(1885.000,796.500){\arc{2850.987}{2.6779}{3.6052}}
\put(-215.000,796.500){\arc{2850.987}{5.8195}{6.7468}}
\put(647.500,571.500){\arc{1277.204}{2.4393}{3.8438}}

\put(722.500,984.000){\arc{1125.000}{3.1416}{5.3559}}


\put(-1124.375,1640.250){\arc{4561.648}{0.2919}{0.7068}}

\put(610,609){\blacken\ellipse{100}{100}}
\put(610,609){\ellipse{100}{100}}
\put(160,159){\blacken\ellipse{100}{100}}
\put(160,159){\ellipse{100}{100}}
\put(610,159){\blacken\ellipse{100}{100}}
\put(610,159){\ellipse{100}{100}}
\put(1060,159){\blacken\ellipse{100}{100}}
\put(1060,159){\ellipse{100}{100}}
\put(160,984){\blacken\ellipse{100}{100}}
\put(160,984){\ellipse{100}{100}}
\put(1060,984){\blacken\ellipse{100}{100}}
\put(1060,984){\ellipse{100}{100}}
\put(1060,1434){\blacken\ellipse{100}{100}}
\put(1060,1434){\ellipse{100}{100}}
\put(610,1434){\blacken\ellipse{100}{100}}
\put(610,1434){\ellipse{100}{100}}

\path(160,169)(610,169)
\path(160,149)(610,149)

\path(610,159)(610,609)
\path(610,159)(1060,159)
\path(160,984)(160,159)

\path(160,994)(1060,994)
\path(160,974)(1060,974)

\path(160,159)(610,609)

\path(1067,167)(617,617)
\path(1053,153)(603,603)

\path(1060,984)(1060,1434)
\path(610,1434)(610,609)

\path(610,1444)(1060,1444)
\path(610,1424)(1060,1424)

\path(610,1434)(1060,984)

\path(610,1434)(160,984)

\put(835,159){\ellipse{70}{70}}
\put(835,159){\whiten\ellipse{70}{70}}
\put(835,159){\blacken\ellipse{30}{30}}
\put(835,159){\ellipse{30}{30}}

\put(1060,1209){\ellipse{70}{70}}
\put(1060,1209){\whiten\ellipse{70}{70}}
\put(1060,1209){\blacken\ellipse{30}{30}}
\put(1060,1209){\ellipse{30}{30}}

\put(387,386){\ellipse{70}{70}}
\put(387,386){\whiten\ellipse{70}{70}}
\put(387,386){\blacken\ellipse{30}{30}}
\put(387,386){\ellipse{30}{30}}

\put(387,1211){\ellipse{70}{70}}
\put(387,1211){\whiten\ellipse{70}{70}}
\put(387,1211){\blacken\ellipse{30}{30}}
\put(387,1211){\ellipse{30}{30}}

\put(20,684){\blacken\ellipse{40}{40}}
\put(20,684){\ellipse{40}{40}}
\put(15,525){\blacken\ellipse{40}{40}}
\put(15,525){\ellipse{40}{40}}

\put(460,800){\blacken\ellipse{40}{40}}
\put(460,800){\ellipse{40}{40}}
\put(470,641){\blacken\ellipse{40}{40}}
\put(470,641){\ellipse{40}{40}}

\put(460,800){\blacken\ellipse{40}{40}}
\put(460,800){\ellipse{40}{40}}
\put(470,641){\blacken\ellipse{40}{40}}
\put(470,641){\ellipse{40}{40}}

\put(1060,1309){\blacken\ellipse{40}{40}}
\put(1060,1309){\ellipse{40}{40}}
\put(1060,1109){\blacken\ellipse{40}{40}}
\put(1060,1109){\ellipse{40}{40}}

\put(760,461){\blacken\ellipse{40}{40}}
\put(760,461){\ellipse{40}{40}}
\put(870,351){\blacken\ellipse{40}{40}}
\put(870,351){\ellipse{40}{40}}

\put(610,459){\whiten\ellipse{50}{50}}
\put(610,459){\ellipse{50}{50}}
\put(610,309){\whiten\ellipse{50}{50}}
\put(610,309){\ellipse{50}{50}}

\put(460,20){\whiten\ellipse{50}{50}}
\put(460,20){\ellipse{50}{50}}
\put(760,20){\whiten\ellipse{50}{50}}
\put(760,20){\ellipse{50}{50}}

\put(920,1130){\whiten\ellipse{50}{50}}
\put(920,1130){\ellipse{50}{50}}
\put(810,1245){\whiten\ellipse{50}{50}}
\put(810,1245){\ellipse{50}{50}}

\put(450,1480){\whiten\ellipse{50}{50}}
\put(450,1480){\ellipse{50}{50}}
\put(610,1540){\whiten\ellipse{50}{50}}
\put(610,1540){\ellipse{50}{50}}

\end{picture}
}\\
$\Omega_3$\\
\end{tabular}
\end{tabular}
$$
which appears in \cite[\S 6.9]{Humphreys90}, although
it is easily
checked that the Gram matrix $G(\Gamma)$ has signature $(1,4)$ and that
the poset $\overline{{\mathscr F}}$ is isomorphic to 
the poset of a combinatorial $4$-simplex.
Thus, by Theorems \ref{Vinberg1}--\ref{Vinberg2},  $\Gamma$ acts cofinitely on $\H^4$ with
fundamental region a cusped $4$-simplex.

From ${\mathscr F}$ and Theorem \ref{Chiswell} we have
$$
\chi(\Gamma)=1/{\mathscr L}(\Gamma),
$$
where ${\mathscr L}(\Gamma)=2^6\,3$, the lowest common multiple of the order of the elements
of $\FF$, is found  using Tables \ref{roots1}--\ref{roots3}.
There are five finite subgroups generated by maximal $S'\subseteq S$: a $D_4$, three
$B_3$'s and an
$A_1\times A_1\times A_1\times A_1$, and
conjugacy class representatives in these groups can be listed using the results
of \S 4. For example,
$x_1,x_2,x_3$ and $x_5$ generate a Weyl group of type $D_4$. Consult the
root system in Table \ref{roots1} and make the
identifications $x_2=s_{\e_1-\e_2}, x_5=s_{\e_2-\e_3},
x_1=s_{\e_3-\e_4}$ and $x_3=s_{\e_3+\e_4}$. By Theorem \ref{Carter2} the order two torsion 
corresponds to the diagrams,
$$
\begin{picture}(10,15)
\put(-12,10){$\e_1-\e_2$}\put(3,3){\circle{6}}\end{picture},
\,\,\,\,\,\,\,\,\,\,\,\,\,\,\,
\begin{picture}(45,15)
\put(-12,10){$\e_1-\e_2$}\put(25,10){$\e_3-\e_4$}
\put(3,3){\circle{6}}
\put(40,3){\circle{6}}\end{picture},
\,\,\,\,\,\,\,\,\,\,\,\,\,\,\,
\begin{picture}(45,15)
\put(-12,10){$\e_1-\e_2$}\put(25,10){$\e_3+\e_4$}
\put(3,3){\circle{6}}
\put(40,3){\circle{6}}\end{picture},
\,\,\,\,\,\,\,\,\,\,\,\,\,\,\,
\begin{picture}(45,15)
\put(-12,10){$\e_3-\e_4$}\put(25,10){$\e_3+\e_4$}
\put(3,3){\circle{6}}
\put(40,3){\circle{6}}\end{picture}
$$
$$
\begin{picture}(84,15)
\put(-12,10){$\e_1-\e_2$}\put(25,10){$\e_3-\e_4$}
\put(63,10){$\e_3+\e_4$}
\put(3,3){\circle{6}}
\put(40,3){\circle{6}}
\put(78,3){\circle{6}}\end{picture},
\,\,\,\,\,\,\text{ and }\,\,\,\,\,\,\,\,\,
\begin{picture}(120,15)
\put(-12,10){$\e_1-\e_2$}\put(25,10){$\e_1+\e_2$}
\put(63,10){$\e_3-\e_4$}\put(101,10){$\e_3+\e_4$}
\put(3,3){\circle{6}}
\put(40,3){\circle{6}}
\put(78,3){\circle{6}}
\put(116,3){\circle{6}}\end{picture},
$$
with resulting representatives $x_2,x_1x_2,x_2x_3,x_1x_3,x_1x_2x_3$ and 
$x_2s_{\e_1+\e_2}x_1x_3$, where $s_{\e_1+\e_2}=wx_3w^{-1}$ with
$w=x_5x_1x_2x_5$; the order three torsion consists of the single class
$$
\begin{picture}(42,15)
\put(-14,10){$\e_1-\e_2$}\put(23,10){$\e_2-\e_3$}
\put(3,3){\circle{6}}\put(6,3){\line(1,0){30}}
\put(39,3){\circle{6}}\end{picture},
$$
with corresponding representative $x_2x_5$. Perform this process for the four other groups.

The three other figures above give actions of the $x_j$ on sets $\om_i$ using the notation from the 
beginning of the section.
By checking that each relator word $(s_{\aa}s_{\bb})^{m_{\aa\bb}}$
acts as the identity on $\Omega_i$, it
can be seen that the diagrams give actions of $\Gamma$ on $\Omega_i$. Let
$U_i$ be the resulting $\C\Gamma$-module, transitive in each case as the diagrams are connected.

One now checks, using Proposition \ref{prop1}, that each representative of torsion
in $\Gamma$ is avoided by at least one of the $U_i$: for instance, those listed 
above are avoided by $U_2,U_2,U_3,U_3,$ $U_3,U_2$ and $U_1$
respectively.
Thus $U=\otimes U_i$ is torsion free. The index
in $\Gamma$ of the corresponding subgroups $\Pi$ is at most
$\prod|\Omega_i\,|=2^6\,3=\LL(\Gamma)$. On the otherhand, the index of
must be a multiple of $\LL(\Gamma)$, hence the
index is $\LL(\Gamma)$. We thus obtain a (non-compact) hyperbolic $4$-manifold
$M=\H^4/\Pi$, with
$\chi(M)=\chi(\Pi)=\LL(\Gamma)\times \chi(\Gamma)=1$. 

\subsection{A $5$-manifold with volume $14\zeta(3)$}

The Coxeter symbol,

$$\begin{picture}(160,45)
\put(40,11){${ x_2}$}\put(2,11){${ x_1}$}\put(82,11){${ x_3}$}
\put(116,11){${ x_4}$}\put(154,11){${ x_5}$}
\put(86,42){${ x_6}$}
\put(23,6){$4$}\put(137,6){$4$}
\put(4,4){\circle{8}}
\put(8,4){\line(1,0){30}}
\put(42,4){\circle{8}}\put(46,4){\line(1,0){30}}
\put(80,4){\circle{8}}\put(84,4){\line(1,0){30}}
\put(118,4){\circle{8}}
\put(122,4){\line(1,0){30}}
\put(156,4){\circle{8}}
\put(80,8){\line(0,1){30}}\put(80,42){\circle{8}}
\end{picture}
$$
corresponds to a $5$-dimensional hyperbolic group acting cofinitely on $\H^5$ with 
fundamental region a cusped $5$-simplex (\cite[\S 6.9]{Humphreys90} or 
Theorems \ref{Vinberg1}--\ref{Vinberg2}). The four diagrams,

$$
\begin{pspicture}(0,0)(14,2.5)
\rput(3,1){\BoxedEPSF{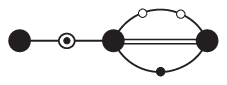 scaled 1000}}
\rput(3,2){\BoxedEPSF{ scaled 1000}}
\rput(6.5,1.5){\BoxedEPSF{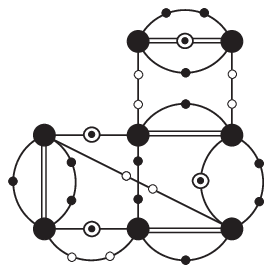 scaled 1000}}
\rput(11,1.5){\BoxedEPSF{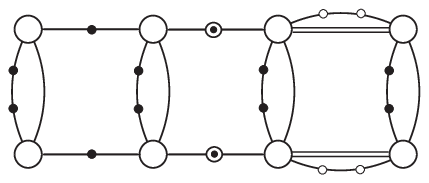 scaled 1000}}
\rput(1.5,2){$\Omega_1$}
\rput(1.5,1){$\Omega_2$}
\rput(5.5,2.5){$\Omega_3$}
\rput(11,.25){$\Omega_4$}
\rput(4.2,1){$\dag$}
\rput(7.7,.4){$\star$}

\rput(9.115,.85){${\sss 1}$}
\rput(10.4,.85){${\sss 3}$}
\rput(11.65,.85){${\sss 5}$}
\rput(12.925,.85){${\sss 7}$}

\rput(9.115,2.125){${\sss 2}$}
\rput(10.375,2.125){${\sss 4}$}
\rput(11.65,2.125){${\sss 6}$}
\rput(12.925,2.125){${\sss 8}$}
\end{pspicture}
$$
\noindent give $\C\Gamma$-modules $U_i$ as in \S 6.1, with the
$\Gamma$-action on $\Omega_4$ imprimitive, each block of imprimitivity
having size eight.
If $x_j$
fixes the block $i$, we use the notation ${\mathbf i}: j\ss,\ldots$ to
mean the action of $x_j$ on the eight points in the block is as a permutation $\ss$ given
by,
$$
\begin{tabular}{l|l|l}
\hline
${\scriptstyle \1:2(2,4)(5,7),3(3,5)(4,6),4(2,5)(4,7)}$&${\scriptstyle \4:2(1,3)(6,8),4(2,4)(5,7)}$
&${\scriptstyle \7:3(3,6)(4,7),6(5,7)(6,8)}$\\
${\scriptstyle \2:2(1,3)(6,8),3(3,4)(5,6),4(2,4)(5,7)}$&${\scriptstyle \5:6(4,6)(5,7)}$
&${\scriptstyle \8:3(1,5)(4,8),6(5,7)(6,8)}$\\
${\scriptstyle \3:2(2,4)(5,7),4(2,5)(4,7)}$&${\scriptstyle \6:6(3,4)(7,8)}$&\\
\hline
\end{tabular}
$$
The permutations $(i,j)_k$ not equal to the identity are,
$$
\begin{tabular}{l|l|l}
\hline
${\scriptstyle (1,2)_1,(3,4)_1=(4,5)}$&${\scriptstyle (7,8)_1,(7,8)_5=(1,3)(5,6)(7,8)}$
&${\scriptstyle (6,8)_2=(2,3,5)(4,7,6)}$\\
${\scriptstyle (1,2)_5,(3,4)_5=(1,2)(3,5,6,4)(7,8)}$&${\scriptstyle (3,5)_3=(3,4,5,6,7,8)}$
&${\scriptstyle (5,7)_2=(2,3)(4,5,6,7,8)}$\\
${\scriptstyle (5,6)_1=(3,6,4)(5,7,8)}$&${\scriptstyle (4,6)_3=(5,7)(6,8)}$
&${\scriptstyle (5,7)_4=(1,2,4,5,6,7,8,3)}$\\
${\scriptstyle (5,6)_5=(1,2)(3,5)(4,7)(6,8)}$&${\scriptstyle (6,8)_4=(1,2,4,8,7,5)(3,6)}$
&${\scriptstyle }$\\
\hline
\end{tabular}
$$
Listing the torsion in $\Gamma$ using the results of \S 4, one can check that the $U_i$ collectively
avoid all torsion, so that $U=\otimes U_i$ is a torsion free $\C\Gamma$-module.

We now compute the volume of the manifold $M^5$ using the results of
\S \ref{modules}. By lemma \ref{lemma3.2}, we have that $U_2\otimes U_3$ is
transitive, so in particular the
size of any $\Gamma$-orbit in $\Omega_1\times \Omega_2\times \Omega_3$
is divisible by $3$. 
On the other hand, using the isomorphism between $\langle x_2,x_3,x_4\rangle$ and the
symmetric group $\SS_4$ given by $x_i\mapsto (i-1,i)$, and
letting
$$
w_1=x_2x_4,w_2=x_3x_2x_3,
$$
we have that $\langle w_1,w_2\rangle$ is a subgroup of order eight. No
non-identity element of it fixes a point of $\Omega_3$, hence of
$\Omega_2\times\Omega_3$. 
However, the subgroup $\langle x_2,x_3,x_4\rangle$ fixes a point
in the first block of imprimitivity of $\Omega_4$.
Proposition \ref{prop2} thus gives that $2^3|\Omega_4|=2^9$
divides the size of any $\Gamma$-orbit in $\Omega_1\times \Omega_2\times
\Omega_3$, hence $2^9\,3=|\Omega_1\times \Omega_2\times\Omega_3|$ does too, and
$U_2\otimes U_3\otimes U_4$ is a transitive $\C\Gamma$-module. 

The subgroup of order two generated by $x_1x_3x_5$ acts on
$\Omega_1$ with no non-identity element fixing a point. On the other
hand, it fixes the point $(\dag,\star,*)$ of $\Omega_2\times
\Omega_3\times \Omega_4$ where $\dag$ and $\star$ are as shown, and $*$
is a point in the seventh block of imprimitivity of $\Omega_4$.
Thus $2|\Omega_2\times
\Omega_3\times \Omega_4|=\times|\Omega_i|$ divides the size of a $\Gamma$-orbit
giving that
$U=\otimes U_i$ is transitive.

There is no Euler characteristic in $5$-dimensions, but by \cite{Ratcliffe99},
a cusped $5$-simplex
fundamental region for $\Gamma$ has volume $7\zeta(3)/2^9\,3$, with $\zeta$ the
Riemann-zeta function. 
Thus $M^5$ has volume $7\zeta(3)\prod|\Omega_i|/2^9\,3=14\zeta(3)$.

\subsection{A $6$-manifold with $|\chi|=16$}

The Coxeter symbol,

$$\begin{picture}(198,45)
\put(40,11){${ x_2}$}\put(2,11){${ x_1}$}\put(82,11){${ x_3}$}
\put(116,11){${ x_4}$}\put(154,11){${ x_5}$}\put(192,11){${ x_6}$}
\put(86,42){${ x_7}$}
\put(175,6){$4$}
\put(4,4){\circle{8}}
\put(8,4){\line(1,0){30}}
\put(42,4){\circle{8}}\put(46,4){\line(1,0){30}}
\put(80,4){\circle{8}}\put(84,4){\line(1,0){30}}
\put(118,4){\circle{8}}
\put(122,4){\line(1,0){30}}\put(156,4){\circle{8}}
\put(160,4){\line(1,0){30}}\put(194,4){\circle{8}}
\put(80,8){\line(0,1){30}}\put(80,42){\circle{8}}
\end{picture}
$$
is, as the others, that of a $\Gamma$ acting cofinitely on $\H^6$ with fundamental region a cusped 
$6$-simplex \cite[\S 6.9]{Humphreys90}. The two diagrams,
$$
\begin{pspicture}(0,0)(14,8)
\rput(1,4){\setlength{\unitlength}{0.00083333in}
\begingroup\makeatletter\ifx\SetFigFont\undefined%
\gdef\SetFigFont#1#2#3#4#5{%
  \reset@font\fontsize{#1}{#2pt}%
  \fontfamily{#3}\fontseries{#4}\fontshape{#5}%
  \selectfont}%
\fi\endgroup%
{\renewcommand{\dashlinestretch}{30}
\begin{picture}(566,483)(0,-10)
\put(283.000,215.250){\arc{487.500}{3.5364}{5.8884}}
\put(283.000,402.750){\arc{487.500}{0.3948}{2.7468}}
\put(283.000,243.375){\arc{468.750}{5.9994}{9.7086}}
\put(58,309){\blacken\ellipse{100}{100}}
\put(58,309){\ellipse{100}{100}}
\put(508,309){\blacken\ellipse{100}{100}}
\put(508,309){\ellipse{100}{100}}
\path(58,319)(508,319)
\path(58,299)(508,299)

\put(283,309){\ellipse{70}{70}}
\put(283,309){\whiten\ellipse{70}{70}}
\put(283,309){\blacken\ellipse{30}{30}}
\put(283,309){\ellipse{30}{30}}

\put(183,30){\blacken\ellipse{40}{40}}
\put(183,30){\ellipse{40}{40}}
\put(383,30){\blacken\ellipse{40}{40}}
\put(383,30){\ellipse{40}{40}}

\put(283,159){\blacken\ellipse{40}{40}}
\put(283,159){\ellipse{40}{40}}

\put(283,9){\whiten\ellipse{50}{50}}
\put(283,9){\ellipse{50}{50}}

\put(183,309){\whiten\ellipse{50}{50}}
\put(183,309){\ellipse{50}{50}}

\put(383,309){\whiten\ellipse{50}{50}}
\put(383,309){\ellipse{50}{50}}

\end{picture}
}}
\rput(9,4){
\setlength{\unitlength}{0.00083333in}
\begingroup\makeatletter\ifx\SetFigFont\undefined%
\gdef\SetFigFont#1#2#3#4#5{%
  \reset@font\fontsize{#1}{#2pt}%
  \fontfamily{#3}\fontseries{#4}\fontshape{#5}%
  \selectfont}%
\fi\endgroup%
{\renewcommand{\dashlinestretch}{30}
\begin{picture}(4811,3909)(0,-10)
\put(1603,2487){\blacken\ellipse{100}{100}}
\put(1603,2487){\ellipse{100}{100}}
\put(1603,2937){\blacken\ellipse{100}{100}}
\put(1603,2937){\ellipse{100}{100}}
\path(1603,2487)(1153,2487)
\path(1603,2937)(1153,2937)
\path(1603,2487)(1603,2937)
\put(1603,2037){\blacken\ellipse{100}{100}}
\put(1603,2037){\ellipse{100}{100}}
\path(1603,2037)(1153,2037)
\path(1603,2037)(1603,2487)
\put(1603,1587){\blacken\ellipse{100}{100}}
\put(1603,1587){\ellipse{100}{100}}
\path(1603,1587)(1153,1587)
\path(1603,1587)(1603,2037)
\put(2953,2037){\blacken\ellipse{100}{100}}
\put(2953,2037){\ellipse{100}{100}}
\put(2953,2487){\blacken\ellipse{100}{100}}
\put(2953,2487){\ellipse{100}{100}}
\path(2953,2037)(2503,2037)
\path(2953,2487)(2503,2487)
\path(2953,2037)(2953,2487)
\put(4753,2037){\blacken\ellipse{100}{100}}
\put(4753,2037){\ellipse{100}{100}}
\put(4753,2487){\blacken\ellipse{100}{100}}
\put(4753,2487){\ellipse{100}{100}}
\path(4753,2037)(4303,2037)
\path(4753,2487)(4303,2487)
\path(4753,2037)(4753,2487)
\put(4753,2937){\blacken\ellipse{100}{100}}
\put(4753,2937){\ellipse{100}{100}}
\put(4303,2937){\blacken\ellipse{100}{100}}
\put(4303,2937){\ellipse{100}{100}}
\path(4753,2937)(4753,2487)
\path(4303,2937)(4303,2487)
\path(4753,2937)(4303,2937)
\put(4753,3387){\blacken\ellipse{100}{100}}
\put(4753,3387){\ellipse{100}{100}}
\put(4303,3387){\blacken\ellipse{100}{100}}
\put(4303,3387){\ellipse{100}{100}}
\path(4753,3387)(4753,2937)
\path(4303,3387)(4303,2937)
\path(4753,3387)(4303,3387)
\put(2503,1587){\blacken\ellipse{100}{100}}
\put(2503,1587){\ellipse{100}{100}}
\path(2503,1587)(2053,1587)
\path(2503,1587)(2503,2037)
\put(2278,1212){\blacken\ellipse{100}{100}}
\put(2278,1212){\ellipse{100}{100}}
\put(2278,1512){\blacken\ellipse{100}{100}}
\put(2278,1512){\ellipse{100}{100}}
\path(2278,1212)(2278,1512)
\put(2278,2112){\blacken\ellipse{100}{100}}
\put(2278,2112){\ellipse{100}{100}}
\put(2278,2412){\blacken\ellipse{100}{100}}
\put(2278,2412){\ellipse{100}{100}}
\path(2278,2112)(2278,2412)
\put(2503,2037){\blacken\ellipse{100}{100}}
\put(2503,2037){\ellipse{100}{100}}
\put(2503,2487){\blacken\ellipse{100}{100}}
\put(2503,2487){\ellipse{100}{100}}
\path(2503,2037)(2053,2037)
\path(2503,2487)(2053,2487)
\path(2503,2037)(2503,2487)

\put(1603,1737){\blacken\ellipse{40}{40}}
\put(1603,1737){\ellipse{40}{40}}
\put(1603,1887){\blacken\ellipse{40}{40}}
\put(1603,1887){\ellipse{40}{40}}

\put(2053,1737){\blacken\ellipse{40}{40}}
\put(2053,1737){\ellipse{40}{40}}
\put(2053,1887){\blacken\ellipse{40}{40}}
\put(2053,1887){\ellipse{40}{40}}

\put(2278,1737){\blacken\ellipse{40}{40}}
\put(2278,1737){\ellipse{40}{40}}
\put(2278,1887){\blacken\ellipse{40}{40}}
\put(2278,1887){\ellipse{40}{40}}

\put(2503,1737){\blacken\ellipse{40}{40}}
\put(2503,1737){\ellipse{40}{40}}
\put(2503,1887){\blacken\ellipse{40}{40}}
\put(2503,1887){\ellipse{40}{40}}

\put(1003,3837){\blacken\ellipse{40}{40}}
\put(1003,3837){\ellipse{40}{40}}

\put(853,3837){\blacken\ellipse{40}{40}}
\put(853,3837){\ellipse{40}{40}}

\put(4753,2637){\blacken\ellipse{40}{40}}
\put(4753,2637){\ellipse{40}{40}}
\put(4753,2787){\blacken\ellipse{40}{40}}
\put(4753,2787){\ellipse{40}{40}}

\put(4303,2637){\blacken\ellipse{40}{40}}
\put(4303,2637){\ellipse{40}{40}}
\put(4303,2787){\blacken\ellipse{40}{40}}
\put(4303,2787){\ellipse{40}{40}}

\put(3853,2637){\blacken\ellipse{40}{40}}
\put(3853,2637){\ellipse{40}{40}}
\put(3853,2787){\blacken\ellipse{40}{40}}
\put(3853,2787){\ellipse{40}{40}}

\put(3403,2637){\blacken\ellipse{40}{40}}
\put(3403,2637){\ellipse{40}{40}}
\put(3403,2787){\blacken\ellipse{40}{40}}
\put(3403,2787){\ellipse{40}{40}}

\put(2953,2637){\blacken\ellipse{40}{40}}
\put(2953,2637){\ellipse{40}{40}}
\put(2953,2787){\blacken\ellipse{40}{40}}
\put(2953,2787){\ellipse{40}{40}}

\path(703,2047)(1153,2047)
\path(703,2027)(1153,2027)
\path(703,2497)(1153,2497)
\path(703,2478)(1153,2478)
\path(703,2947)(1153,2947)
\path(703,2927)(1153,2927)

\path(703,1597)(1153,1597)
\path(703,1577)(1153,1577)

\path(703,1147)(1153,1147)
\path(703,1127)(1153,1127)

\path(703,697)(1153,697)
\path(703,677)(1153,677)

\path(703,247)(1153,247)
\path(703,227)(1153,227)

\put(-740.750,912.000){\arc{3187.500}{5.8458}{6.7205}}
\put(3046.750,912.000){\arc{3187.500}{2.7043}{3.5789}}
\put(2728.000,1044.000){\arc{3186.141}{4.2749}{5.1499}}
\put(2728.000,3480.000){\arc{3186.141}{1.1333}{2.0083}}
\put(2973.833,4374.500){\arc{4164.416}{1.1349}{1.9115}}
\put(2973.833,149.500){\arc{4164.416}{4.3717}{5.1483}}
\put(-290.750,912.000){\arc{3187.500}{5.8458}{6.7205}}
\put(1003.000,1137.000){\arc{1500.000}{5.6397}{6.9267}}
\put(1003.000,687.000){\arc{1500.000}{5.6397}{6.9267}}
\put(1753.000,687.000){\arc{1500.000}{2.4981}{3.7851}}
\put(1753.000,1137.000){\arc{1500.000}{2.4981}{3.7851}}
\put(1303.000,1137.000){\arc{1500.000}{2.4981}{3.7851}}
\put(1303.000,687.000){\arc{1500.000}{2.4981}{3.7851}}
\put(553.000,2299.500){\arc{604.669}{1.0517}{3.0172}}
\put(2953.000,2262.000){\arc{5850.000}{2.7468}{3.5364}}
\put(215.500,2374.500){\arc{2186.607}{5.7428}{6.3863}}
\put(2540.500,2374.500){\arc{2186.607}{3.0385}{3.6820}}
\put(1828.000,143.250){\arc{3187.500}{4.2751}{5.1497}}
\put(1828.000,3480.750){\arc{3187.500}{1.1335}{2.0081}}
\put(1828.000,1043.250){\arc{3187.500}{4.2751}{5.1497}}
\put(-459.500,1924.500){\arc{3589.046}{6.0940}{6.7375}}
\put(3215.500,1924.500){\arc{3589.046}{2.6873}{3.3308}}
\put(1331.250,2093.250){\arc{2238.093}{2.1170}{3.2930}}

\put(253,3837){\blacken\ellipse{100}{100}}
\put(253,3837){\ellipse{100}{100}}

\put(703,3387){\blacken\ellipse{100}{100}}
\put(703,3387){\ellipse{100}{100}}
\put(703,2937){\blacken\ellipse{100}{100}}
\put(703,2937){\ellipse{100}{100}}

\put(1153,2937){\blacken\ellipse{100}{100}}
\put(1153,2937){\ellipse{100}{100}}
\put(1153,3837){\blacken\ellipse{100}{100}}
\put(1153,3837){\ellipse{100}{100}}

\put(1603,3837){\blacken\ellipse{100}{100}}
\put(1603,3837){\ellipse{100}{100}}
\put(2503,1137){\blacken\ellipse{100}{100}}
\put(2503,1137){\ellipse{100}{100}}
\put(2053,1137){\blacken\ellipse{100}{100}}
\put(2053,1137){\ellipse{100}{100}}

\put(478,3837){\blacken\ellipse{40}{40}}
\put(478,3837){\ellipse{40}{40}}
\put(478,3387){\blacken\ellipse{40}{40}}
\put(478,3387){\ellipse{40}{40}}

\put(1378,3837){\blacken\ellipse{40}{40}}
\put(1378,3837){\ellipse{40}{40}}

\put(2053,2262){\blacken\ellipse{40}{40}}
\put(2053,2262){\ellipse{40}{40}}

\put(1303,2262){\blacken\ellipse{100}{100}}
\put(1303,2262){\ellipse{100}{100}}

\put(1453,2262){\blacken\ellipse{100}{100}}
\put(1453,2262){\ellipse{100}{100}}

\put(1603,2262){\blacken\ellipse{40}{40}}
\put(1603,2262){\ellipse{40}{40}}

\put(2053,2037){\blacken\ellipse{100}{100}}
\put(2053,2037){\ellipse{100}{100}}
\put(1603,912){\blacken\ellipse{40}{40}}
\put(1603,912){\ellipse{40}{40}}
\put(1453,912){\blacken\ellipse{40}{40}}
\put(1453,912){\ellipse{40}{40}}

\put(1303,912){\blacken\ellipse{40}{40}}
\put(1303,912){\ellipse{40}{40}}

\put(1153,912){\blacken\ellipse{40}{40}}
\put(1153,912){\ellipse{40}{40}}
\put(853,912){\blacken\ellipse{40}{40}}
\put(853,912){\ellipse{40}{40}}
\put(703,912){\blacken\ellipse{40}{40}}
\put(703,912){\ellipse{40}{40}}
\put(253,912){\blacken\ellipse{40}{40}}
\put(253,912){\ellipse{40}{40}}

\put(2053,1362){\blacken\ellipse{40}{40}}
\put(2053,1362){\ellipse{40}{40}}
\put(2278,1362){\blacken\ellipse{40}{40}}
\put(2278,1362){\ellipse{40}{40}}
\put(2503,2262){\blacken\ellipse{40}{40}}
\put(2503,2262){\ellipse{40}{40}}
\put(2953,2262){\blacken\ellipse{40}{40}}
\put(2953,2262){\ellipse{40}{40}}
\put(3403,2262){\blacken\ellipse{40}{40}}
\put(3403,2262){\ellipse{40}{40}}
\put(3853,2262){\blacken\ellipse{40}{40}}
\put(3853,2262){\ellipse{40}{40}}
\put(4303,2262){\blacken\ellipse{40}{40}}
\put(4303,2262){\ellipse{40}{40}}
\put(4753,2262){\blacken\ellipse{40}{40}}
\put(4753,2262){\ellipse{40}{40}}
\put(3403,3162){\blacken\ellipse{40}{40}}
\put(3403,3162){\ellipse{40}{40}}
\put(2953,3162){\blacken\ellipse{40}{40}}
\put(2953,3162){\ellipse{40}{40}}
\put(3853,3162){\blacken\ellipse{40}{40}}
\put(3853,3162){\ellipse{40}{40}}
\put(4303,3162){\blacken\ellipse{40}{40}}
\put(4303,3162){\ellipse{40}{40}}
\put(4753,3162){\blacken\ellipse{40}{40}}
\put(4753,3162){\ellipse{40}{40}}
\put(2503,1362){\blacken\ellipse{40}{40}}
\put(2503,1362){\ellipse{40}{40}}
\put(2053,1587){\blacken\ellipse{100}{100}}
\put(2053,1587){\ellipse{100}{100}}

\put(3403,3387){\blacken\ellipse{100}{100}}
\put(3403,3387){\ellipse{100}{100}}
\put(2953,3387){\blacken\ellipse{100}{100}}
\put(2953,3387){\ellipse{100}{100}}
\put(3853,3387){\blacken\ellipse{100}{100}}
\put(3853,3387){\ellipse{100}{100}}
\put(2953,2937){\blacken\ellipse{100}{100}}
\put(2953,2937){\ellipse{100}{100}}
\put(3403,2937){\blacken\ellipse{100}{100}}
\put(3403,2937){\ellipse{100}{100}}
\put(3853,2937){\blacken\ellipse{100}{100}}
\put(3853,2937){\ellipse{100}{100}}
\put(3853,2037){\blacken\ellipse{100}{100}}
\put(3853,2037){\ellipse{100}{100}}
\put(3853,2487){\blacken\ellipse{100}{100}}
\put(3853,2487){\ellipse{100}{100}}
\put(3403,2037){\blacken\ellipse{100}{100}}
\put(3403,2037){\ellipse{100}{100}}
\put(3403,2487){\blacken\ellipse{100}{100}}
\put(3403,2487){\ellipse{100}{100}}
\put(4303,2037){\blacken\ellipse{100}{100}}
\put(4303,2037){\ellipse{100}{100}}
\put(4303,2487){\blacken\ellipse{100}{100}}
\put(4303,2487){\ellipse{100}{100}}

\put(253,3387){\blacken\ellipse{100}{100}}
\put(253,3387){\ellipse{100}{100}}
\put(703,3837){\blacken\ellipse{100}{100}}
\put(703,3837){\ellipse{100}{100}}

\put(28,2337){\blacken\ellipse{40}{40}}
\put(28,2337){\ellipse{40}{40}}
\put(28,2187){\blacken\ellipse{40}{40}}
\put(28,2187){\ellipse{40}{40}}

\put(1255,2712){\blacken\ellipse{40}{40}}
\put(1255,2712){\ellipse{40}{40}}

\put(1330,1887){\blacken\ellipse{40}{40}}
\put(1330,1887){\ellipse{40}{40}}
\put(1320,1737){\blacken\ellipse{40}{40}}
\put(1320,1737){\ellipse{40}{40}}

\put(1501,2712){\blacken\ellipse{40}{40}}
\put(1501,2712){\ellipse{40}{40}}
\put(1430,1887){\blacken\ellipse{40}{40}}
\put(1430,1887){\ellipse{40}{40}}
\put(1440,1737){\blacken\ellipse{40}{40}}
\put(1440,1737){\ellipse{40}{40}}

\put(703,237){\blacken\ellipse{100}{100}}
\put(703,237){\ellipse{100}{100}}
\put(1153,237){\blacken\ellipse{100}{100}}
\put(1153,237){\ellipse{100}{100}}
\put(1603,237){\blacken\ellipse{100}{100}}
\put(1603,237){\ellipse{100}{100}}
\put(1603,687){\blacken\ellipse{100}{100}}
\put(1603,687){\ellipse{100}{100}}
\put(1603,1137){\blacken\ellipse{100}{100}}
\put(1603,1137){\ellipse{100}{100}}
\put(253,1137){\blacken\ellipse{100}{100}}
\put(253,1137){\ellipse{100}{100}}
\put(253,687){\blacken\ellipse{100}{100}}
\put(253,687){\ellipse{100}{100}}
\put(703,1137){\blacken\ellipse{100}{100}}
\put(703,1137){\ellipse{100}{100}}
\put(1153,1137){\blacken\ellipse{100}{100}}
\put(1153,1137){\ellipse{100}{100}}

\put(703,687){\blacken\ellipse{100}{100}}
\put(703,687){\ellipse{100}{100}}
\put(1153,687){\blacken\ellipse{100}{100}}
\put(1153,687){\ellipse{100}{100}}

\put(703,1737){\blacken\ellipse{40}{40}}
\put(703,1737){\ellipse{40}{40}}
\put(703,1887){\blacken\ellipse{40}{40}}
\put(703,1887){\ellipse{40}{40}}

\put(703,1587){\blacken\ellipse{100}{100}}
\put(703,1587){\ellipse{100}{100}}
\put(1153,1587){\blacken\ellipse{100}{100}}
\put(1153,1587){\ellipse{100}{100}}

\put(1153,1737){\blacken\ellipse{40}{40}}
\put(1153,1737){\ellipse{40}{40}}
\put(1153,1887){\blacken\ellipse{40}{40}}
\put(1153,1887){\ellipse{40}{40}}

\put(703,2037){\blacken\ellipse{100}{100}}
\put(703,2037){\ellipse{100}{100}}
\put(1153,2037){\blacken\ellipse{100}{100}}
\put(1153,2037){\ellipse{100}{100}}
\put(703,2487){\blacken\ellipse{100}{100}}
\put(703,2487){\ellipse{100}{100}}
\put(1153,2487){\blacken\ellipse{100}{100}}
\put(1153,2487){\ellipse{100}{100}}
\put(2053,2487){\blacken\ellipse{100}{100}}
\put(2053,2487){\ellipse{100}{100}}
\put(2278,2262){\blacken\ellipse{40}{40}}
\put(2278,2262){\ellipse{40}{40}}

\put(253,2262){\blacken\ellipse{100}{100}}
\put(253,2262){\ellipse{100}{100}}

\put(478,2600){\blacken\ellipse{40}{40}}
\put(478,2600){\ellipse{40}{40}}

\put(253,1812){\blacken\ellipse{40}{40}}
\put(253,1812){\ellipse{40}{40}}

\put(300,1662){\blacken\ellipse{40}{40}}
\put(300,1662){\ellipse{40}{40}}

\put(703,2337){\blacken\ellipse{40}{40}}
\put(703,2337){\ellipse{40}{40}}
\put(1153,2337){\blacken\ellipse{40}{40}}
\put(1153,2337){\ellipse{40}{40}}

\path(253,3837)(703,3837)
\path(703,2937)(703,3387)
\path(4303,3387)(4753,2937)
\path(4753,3387)(4303,2937)
\path(2503,1137)(2503,1587)
\path(2503,1137)(2053,1137)
\path(2053,1137)(2053,1587)

\path(2053,1137)(2278,1212)
\path(2053,1157)(2278,1232)
\path(2053,1567)(2278,1492)
\path(2053,1587)(2278,1512)
\path(2053,2037)(2278,2112)
\path(2053,2057)(2278,2132)
\path(2053,2467)(2278,2392)
\path(2053,2487)(2278,2412)

\path(2278,2112)(2278,1512)
\path(1603,237)(1603,87)(2053,87)(2053,1137)
\path(1153,237)(1153,12)(2503,12)(2503,1137)
\path(4265,3425)(4715,2975)
\path(253,3387)(703,3387)

\path(1153,3837)(1603,3837)

\path(2053,2037)(2053,2487)
\path(1153,2037)(1303,2262)
\path(1453,2262)(1603,2037)
\path(1303,2262)(1453,2262)

\path(703,3837)(1153,3837)

\path(3403,3387)(3403,2937)
\path(2953,3387)(2953,2937)
\path(3403,3387)(2953,3387)
\path(3853,3387)(3853,2937)
\path(4303,3387)(3853,3387)
\path(3403,2937)(2953,2937)
\path(3403,2937)(3403,2487)
\path(3853,2937)(3853,2487)
\path(4303,2937)(3853,2937)

\path(3853,2037)(3853,2487)
\path(3403,2037)(2953,2037)
\path(3403,2487)(2953,2487)
\path(3403,2037)(3403,2487)
\path(4303,2037)(3853,2037)
\path(4303,2487)(3853,2487)
\path(4303,2037)(4303,2487)

\path(3853,3397)(3403,3397)
\path(3853,3378)(3403,3378)

\path(3853,2947)(3403,2947)
\path(3853,2927)(3403,2927)

\path(3853,2497)(3403,2497)
\path(3853,2478)(3403,2478)

\path(3853,2047)(3403,2047)
\path(3853,2027)(3403,2027)

\path(2053,1587)(1603,1587)
\path(2053,1587)(2053,2037)
\path(253,3837)(253,3387)
\path(703,3837)(703,3387)
\path(2953,2937)(2953,2487)
\path(253,3387)(253,2262)

\path(703,237)(703,687)
\path(1153,237)(1153,687)
\path(703,687)(703,1137)
\path(1153,687)(1153,1137)
\path(1603,237)(1153,237)
\path(1603,237)(1603,687)
\path(1603,687)(1153,687)
\path(1603,687)(1603,1137)
\path(1603,1137)(1153,1137)
\path(1603,1137)(1603,1587)
\path(253,1137)(703,1137)
\path(253,687)(703,687)
\path(253,1137)(253,687)
\path(703,1137)(703,1587)
\path(1153,1137)(1153,1587)
\path(703,1587)(703,2037)
\path(1153,1587)(1153,2037)
\path(703,2037)(703,2487)
\path(1153,2037)(1153,2487)
\path(703,2487)(703,2937)
\path(1153,2487)(1153,2937)
\path(2053,2037)(1603,2037)
\path(2053,2487)(1603,2487)

\put(478,1140){\ellipse{70}{70}}
\put(478,1140){\whiten\ellipse{70}{70}}
\put(478,1140){\blacken\ellipse{30}{30}}
\put(478,1140){\ellipse{30}{30}}

\put(478,690){\ellipse{70}{70}}
\put(478,690){\whiten\ellipse{70}{70}}
\put(478,690){\blacken\ellipse{30}{30}}
\put(478,690){\ellipse{30}{30}}

\put(251,2862){\ellipse{70}{70}}
\put(251,2862){\whiten\ellipse{70}{70}}
\put(251,2862){\blacken\ellipse{30}{30}}
\put(251,2862){\ellipse{30}{30}}

\put(701,3162){\ellipse{70}{70}}
\put(701,3162){\whiten\ellipse{70}{70}}
\put(701,3162){\blacken\ellipse{30}{30}}
\put(701,3162){\ellipse{30}{30}}

\put(1828,2490){\ellipse{70}{70}}
\put(1828,2490){\whiten\ellipse{70}{70}}
\put(1828,2490){\blacken\ellipse{30}{30}}
\put(1828,2490){\ellipse{30}{30}}

\put(1828,2040){\ellipse{70}{70}}
\put(1828,2040){\whiten\ellipse{70}{70}}
\put(1828,2040){\blacken\ellipse{30}{30}}
\put(1828,2040){\ellipse{30}{30}}

\put(1828,1590){\ellipse{70}{70}}
\put(1828,1590){\whiten\ellipse{70}{70}}
\put(1828,1590){\blacken\ellipse{30}{30}}
\put(1828,1590){\ellipse{30}{30}}

\put(1828,1740){\ellipse{70}{70}}
\put(1828,1740){\whiten\ellipse{70}{70}}
\put(1828,1740){\blacken\ellipse{30}{30}}
\put(1828,1740){\ellipse{30}{30}}

\put(1828,1890){\ellipse{70}{70}}
\put(1828,1890){\whiten\ellipse{70}{70}}
\put(1828,1890){\blacken\ellipse{30}{30}}
\put(1828,1890){\ellipse{30}{30}}

\put(1828,2640){\ellipse{70}{70}}
\put(1828,2640){\whiten\ellipse{70}{70}}
\put(1828,2640){\blacken\ellipse{30}{30}}
\put(1828,2640){\ellipse{30}{30}}

\put(4078,2039){\ellipse{70}{70}}
\put(4078,2039){\whiten\ellipse{70}{70}}
\put(4078,2039){\blacken\ellipse{30}{30}}
\put(4078,2039){\ellipse{30}{30}}

\put(4078,2489){\ellipse{70}{70}}
\put(4078,2489){\whiten\ellipse{70}{70}}
\put(4078,2489){\blacken\ellipse{30}{30}}
\put(4078,2489){\ellipse{30}{30}}

\put(4078,2939){\ellipse{70}{70}}
\put(4078,2939){\whiten\ellipse{70}{70}}
\put(4078,2939){\blacken\ellipse{30}{30}}
\put(4078,2939){\ellipse{30}{30}}

\put(4078,3389){\ellipse{70}{70}}
\put(4078,3389){\whiten\ellipse{70}{70}}
\put(4078,3389){\blacken\ellipse{30}{30}}
\put(4078,3389){\ellipse{30}{30}}

\put(2056,687){\ellipse{70}{70}}
\put(2056,687){\whiten\ellipse{70}{70}}
\put(2056,687){\blacken\ellipse{30}{30}}
\put(2056,687){\ellipse{30}{30}}

\put(2506,687){\ellipse{70}{70}}
\put(2506,687){\whiten\ellipse{70}{70}}
\put(2506,687){\blacken\ellipse{30}{30}}
\put(2506,687){\ellipse{30}{30}}


\path(703,2937)(253,2262)

\path(253,2272)(1303,2272)
\path(253,2252)(1303,2252)

\put(2653,2487){\whiten\ellipse{50}{50}}
\put(2653,2487){\ellipse{50}{50}}
\put(703,3687){\whiten\ellipse{50}{50}}
\put(703,3687){\ellipse{50}{50}}
\put(703,3537){\whiten\ellipse{50}{50}}
\put(703,3537){\ellipse{50}{50}}
\put(253,3687){\whiten\ellipse{50}{50}}
\put(253,3687){\ellipse{50}{50}}
\put(253,3537){\whiten\ellipse{50}{50}}
\put(253,3537){\ellipse{50}{50}}
\put(1153,2712){\whiten\ellipse{50}{50}}
\put(1153,2712){\ellipse{50}{50}}
\put(4528,3387){\whiten\ellipse{50}{50}}
\put(4528,3387){\ellipse{50}{50}}
\put(4528,2937){\whiten\ellipse{50}{50}}
\put(4528,2937){\ellipse{50}{50}}
\put(4528,2487){\whiten\ellipse{50}{50}}
\put(4528,2487){\ellipse{50}{50}}
\put(4528,2037){\whiten\ellipse{50}{50}}
\put(4528,2037){\ellipse{50}{50}}
\put(703,462){\whiten\ellipse{50}{50}}
\put(703,462){\ellipse{50}{50}}
\put(1153,462){\whiten\ellipse{50}{50}}
\put(1153,462){\ellipse{50}{50}}
\put(1603,462){\whiten\ellipse{50}{50}}
\put(1603,462){\ellipse{50}{50}}
\put(703,1362){\whiten\ellipse{50}{50}}
\put(703,1362){\ellipse{50}{50}}
\put(1153,1362){\whiten\ellipse{50}{50}}
\put(1153,1362){\ellipse{50}{50}}
\put(1603,1362){\whiten\ellipse{50}{50}}
\put(1603,1362){\ellipse{50}{50}}
\put(553,762){\whiten\ellipse{50}{50}}
\put(553,762){\ellipse{50}{50}}
\put(553,612){\whiten\ellipse{50}{50}}
\put(553,612){\ellipse{50}{50}}
\put(553,1062){\whiten\ellipse{50}{50}}
\put(553,1062){\ellipse{50}{50}}
\put(553,1212){\whiten\ellipse{50}{50}}
\put(553,1212){\ellipse{50}{50}}
\put(1003,612){\whiten\ellipse{50}{50}}
\put(1003,612){\ellipse{50}{50}}
\put(1753,612){\whiten\ellipse{50}{50}}
\put(1753,612){\ellipse{50}{50}}
\put(1753,762){\whiten\ellipse{50}{50}}
\put(1753,762){\ellipse{50}{50}}
\put(1003,1212){\whiten\ellipse{50}{50}}
\put(1003,1212){\ellipse{50}{50}}
\put(1003,1062){\whiten\ellipse{50}{50}}
\put(1003,1062){\ellipse{50}{50}}
\put(1753,1212){\whiten\ellipse{50}{50}}
\put(1753,1212){\ellipse{50}{50}}
\put(1753,1062){\whiten\ellipse{50}{50}}
\put(1753,1062){\ellipse{50}{50}}
\put(703,2712){\whiten\ellipse{50}{50}}
\put(703,2712){\ellipse{50}{50}}
\put(1603,2712){\whiten\ellipse{50}{50}}
\put(1603,2712){\ellipse{50}{50}}
\put(4453,3237){\whiten\ellipse{50}{50}}
\put(4453,3237){\ellipse{50}{50}}
\put(4603,3087){\whiten\ellipse{50}{50}}
\put(4603,3087){\ellipse{50}{50}}
\put(4453,3087){\whiten\ellipse{50}{50}}
\put(4453,3087){\ellipse{50}{50}}
\put(4603,3237){\whiten\ellipse{50}{50}}
\put(4603,3237){\ellipse{50}{50}}
\put(2653,2637){\whiten\ellipse{50}{50}}
\put(2653,2637){\ellipse{50}{50}}
\put(2803,2637){\whiten\ellipse{50}{50}}
\put(2803,2637){\ellipse{50}{50}}
\put(2803,1887){\whiten\ellipse{50}{50}}
\put(2803,1887){\ellipse{50}{50}}
\put(2803,2300){\whiten\ellipse{50}{50}}
\put(2803,2300){\ellipse{50}{50}}
\put(2803,2220){\whiten\ellipse{50}{50}}
\put(2803,2220){\ellipse{50}{50}}
\put(2803,2487){\whiten\ellipse{50}{50}}
\put(2803,2487){\ellipse{50}{50}}
\put(2803,2037){\whiten\ellipse{50}{50}}
\put(2803,2037){\ellipse{50}{50}}
\put(2653,1887){\whiten\ellipse{50}{50}}
\put(2653,1887){\ellipse{50}{50}}
\put(2653,2037){\whiten\ellipse{50}{50}}
\put(2653,2037){\ellipse{50}{50}}
\put(3103,2298){\whiten\ellipse{50}{50}}
\put(3103,2298){\ellipse{50}{50}}
\put(3103,2220){\whiten\ellipse{50}{50}}
\put(3103,2220){\ellipse{50}{50}}
\put(403,2037){\whiten\ellipse{50}{50}}
\put(403,2037){\ellipse{50}{50}}
\put(1003,762){\whiten\ellipse{50}{50}}
\put(1003,762){\ellipse{50}{50}}

\put(1520,2160){\whiten\ellipse{50}{50}}
\put(1520,2160){\ellipse{50}{50}}

\put(1235,2160){\whiten\ellipse{50}{50}}
\put(1235,2160){\ellipse{50}{50}}
\end{picture}
}
}
\rput(1,3){$\Omega_1$}\rput(8,7){$\Omega_2$}
\rput(7.6,8){$\dag$}
\end{pspicture}
$$
\noindent give primitive $\Gamma$-actions on $2$ and $64$ points. Diagram $\om_3$ below depicts an action 
with $27$ blocks of imprimitivity, each of size 
six\footnote{Nikolai Vavilov has pointed out to me the curious fact that $\Omega_3$ is the weight diagram of
the representation of the Chevalley group of type $E_6$ with highest weight (see for example \cite{Vavilov}).}:
$$
\begin{pspicture}(0,0)(14,3.5)
\rput(7,0.3){$\Omega_3$}
\rput(7,1.75){
\setlength{\unitlength}{0.00083333in}
\begingroup\makeatletter\ifx\SetFigFont\undefined%
\gdef\SetFigFont#1#2#3#4#5{%
  \reset@font\fontsize{#1}{#2pt}%
  \fontfamily{#3}\fontseries{#4}\fontshape{#5}%
  \selectfont}%
\fi\endgroup%
{\renewcommand{\dashlinestretch}{30}
\begin{picture}(5446,1663)(0,-10)

\put(2800,1080){\blacken\ellipse{40}{40}}
\put(2800,1080){\ellipse{40}{40}}
\put(2900,1180){\blacken\ellipse{40}{40}}
\put(2900,1180){\ellipse{40}{40}}
\put(3100,780){\blacken\ellipse{40}{40}}
\put(3100,780){\ellipse{40}{40}}
\put(3200,880){\blacken\ellipse{40}{40}}
\put(3200,880){\ellipse{40}{40}}
\put(3400,480){\blacken\ellipse{40}{40}}
\put(3400,480){\ellipse{40}{40}}
\put(3500,580){\blacken\ellipse{40}{40}}
\put(3500,580){\ellipse{40}{40}}
\put(3700,180){\blacken\ellipse{40}{40}}
\put(3700,180){\ellipse{40}{40}}
\put(3800,280){\blacken\ellipse{40}{40}}
\put(3800,280){\ellipse{40}{40}}

\put(2500,1380){\blacken\ellipse{40}{40}}
\put(2500,1380){\ellipse{40}{40}}
\put(2600,1480){\blacken\ellipse{40}{40}}
\put(2600,1480){\ellipse{40}{40}}

\put(200,673){\blacken\ellipse{40}{40}}
\put(200,673){\ellipse{40}{40}}
\put(320,673){\blacken\ellipse{40}{40}}
\put(320,673){\ellipse{40}{40}}


\put(73,673){\whiten\ellipse{130}{130}}
\put(73,673){\ellipse{130}{130}}
\put(448,673){\whiten\ellipse{130}{130}}
\put(448,673){\ellipse{130}{130}}
\put(1198,673){\whiten\ellipse{130}{130}}
\put(1198,673){\ellipse{130}{130}}
\put(1498,973){\whiten\ellipse{130}{130}}
\put(1498,973){\ellipse{130}{130}}
\put(1498,373){\whiten\ellipse{130}{130}}
\put(1498,373){\ellipse{130}{130}}
\put(1798,73){\whiten\ellipse{130}{130}}
\put(1798,73){\ellipse{130}{130}}
\put(2098,373){\whiten\ellipse{130}{130}}
\put(2098,373){\ellipse{130}{130}}
\put(2098,973){\whiten\ellipse{130}{130}}
\put(2098,973){\ellipse{130}{130}}
\put(2398,673){\whiten\ellipse{130}{130}}
\put(2398,673){\ellipse{130}{130}}
\put(2698,973){\whiten\ellipse{130}{130}}
\put(2698,973){\ellipse{130}{130}}
\put(2398,1273){\whiten\ellipse{130}{130}}
\put(2398,1273){\ellipse{130}{130}}
\put(2698,1573){\whiten\ellipse{130}{130}}
\put(2698,1573){\ellipse{130}{130}}
\put(2998,1273){\whiten\ellipse{130}{130}}
\put(2998,1273){\ellipse{130}{130}}
\put(3298,973){\whiten\ellipse{130}{130}}
\put(3298,973){\ellipse{130}{130}}
\put(3298,373){\whiten\ellipse{130}{130}}
\put(3298,373){\ellipse{130}{130}}
\put(3598,73){\whiten\ellipse{130}{130}}
\put(3598,73){\ellipse{130}{130}}
\put(3898,973){\whiten\ellipse{130}{130}}
\put(3898,973){\ellipse{130}{130}}
\put(4198,673){\whiten\ellipse{130}{130}}
\put(4198,673){\ellipse{130}{130}}
\put(4948,673){\whiten\ellipse{130}{130}}
\put(4948,673){\ellipse{130}{130}}
\put(5323,673){\whiten\ellipse{130}{130}}
\put(5323,673){\ellipse{130}{130}}
\put(823,673){\whiten\ellipse{130}{130}}
\put(823,673){\ellipse{130}{130}}
\put(4573,673){\whiten\ellipse{130}{130}}
\put(4573,673){\ellipse{130}{130}}
\put(3898,373){\whiten\ellipse{130}{130}}
\put(3898,373){\ellipse{130}{130}}
\put(3598,673){\whiten\ellipse{130}{130}}
\put(3598,673){\ellipse{130}{130}}
\put(2698,373){\whiten\ellipse{130}{130}}
\put(2698,373){\ellipse{130}{130}}
\put(2998,673){\whiten\ellipse{130}{130}}
\put(2998,673){\ellipse{130}{130}}
\put(1798,673){\whiten\ellipse{130}{130}}
\put(1798,673){\ellipse{130}{130}}

\path(145,673)(373,673)
\path(520,673)(751,673)
\path(895,673)(1125,673)

\path(4270,673)(4498,673)

\path(4648,685)(4873,685)
\path(4648,661)(4873,661)

\path(5023,673)(5248,673)

\path(3942,418)(4155,630)
\path(3642,118)(3855,330)
\path(3342,418)(3555,630)
\path(3642,718)(3855,930)
\path(2742,418)(2955,630)
\path(3042,718)(3255,930)
\path(2442,718)(2655,930)
\path(2142,1018)(2355,1230)
\path(2742,1018)(2955,1230)
\path(2442,1318)(2655,1530)
\path(1842,718)(2055,930)
\path(2142,418)(2355,630)
\path(1542,418)(1755,630)
\path(1842,118)(2055,330)
\path(1242,718)(1455,930)
\path(4155,718)(3935,925)
\path(3855,418)(3635,625)
\path(3555,118)(3335,325)
\path(3555,718)(3335,925)

\path(3265,1028)(3045,1235)
\path(3245,1008)(3025,1215)

\path(3255,418)(3035,625)

\path(2965,728)(2745,935)
\path(2945,708)(2725,915)

\path(2665,428)(2445,635)
\path(2645,408)(2425,615)

\path(2655,1018)(2435,1225)
\path(2955,1318)(2735,1525)
\path(2355,718)(2135,925)
\path(2055,418)(1835,625)
\path(1755,118)(1535,325)

\path(1465,428)(1245,635)
\path(1445,408)(1225,615)

\path(1765,728)(1545,935)
\path(1745,708)(1525,915)

\put(50,650){\makebox(0,0)[lb]{\smash{{{\SetFigFont{5}{6.0}{\rmdefault}{\mddefault}{\updefault}$1$}}}}}
\put(425,650){\makebox(0,0)[lb]{\smash{{{\SetFigFont{5}{6.0}{\rmdefault}{\mddefault}{\updefault}$2$}}}}}
\put(800,650){\makebox(0,0)[lb]{\smash{{{\SetFigFont{5}{6.0}{\rmdefault}{\mddefault}{\updefault}$3$}}}}}
\put(1175,650){\makebox(0,0)[lb]{\smash{{{\SetFigFont{5}{6.0}{\rmdefault}{\mddefault}{\updefault}$4$}}}}}
\put(1475,950){\makebox(0,0)[lb]{\smash{{{\SetFigFont{5}{6.0}{\rmdefault}{\mddefault}{\updefault}$5$}}}}}
\put(1475,350){\makebox(0,0)[lb]{\smash{{{\SetFigFont{5}{6.0}{\rmdefault}{\mddefault}{\updefault}$6$}}}}}
\put(1775,650){\makebox(0,0)[lb]{\smash{{{\SetFigFont{5}{6.0}{\rmdefault}{\mddefault}{\updefault}$7$}}}}}
\put(1775,50){\makebox(0,0)[lb]{\smash{{{\SetFigFont{5}{6.0}{\rmdefault}{\mddefault}{\updefault}$8$}}}}}
\put(2075,950){\makebox(0,0)[lb]{\smash{{{\SetFigFont{5}{6.0}{\rmdefault}{\mddefault}{\updefault}$9$}}}}}
\put(2045,350){\makebox(0,0)[lb]{\smash{{{\SetFigFont{5}{6.0}{\rmdefault}{\mddefault}{\updefault}$10$}}}}}
\put(2345,650){\makebox(0,0)[lb]{\smash{{{\SetFigFont{5}{6.0}{\rmdefault}{\mddefault}{\updefault}$12$}}}}}
\put(2345,1250){\makebox(0,0)[lb]{\smash{{{\SetFigFont{5}{6.0}{\rmdefault}{\mddefault}{\updefault}$11$}}}}}
\put(2645,1550){\makebox(0,0)[lb]{\smash{{{\SetFigFont{5}{6.0}{\rmdefault}{\mddefault}{\updefault}$13$}}}}}
\put(2645,950){\makebox(0,0)[lb]{\smash{{{\SetFigFont{5}{6.0}{\rmdefault}{\mddefault}{\updefault}$14$}}}}}
\put(2645,350){\makebox(0,0)[lb]{\smash{{{\SetFigFont{5}{6.0}{\rmdefault}{\mddefault}{\updefault}$15$}}}}}
\put(2945,1250){\makebox(0,0)[lb]{\smash{{{\SetFigFont{5}{6.0}{\rmdefault}{\mddefault}{\updefault}$16$}}}}}
\put(2945,650){\makebox(0,0)[lb]{\smash{{{\SetFigFont{5}{6.0}{\rmdefault}{\mddefault}{\updefault}$17$}}}}}
\put(3245,950){\makebox(0,0)[lb]{\smash{{{\SetFigFont{5}{6.0}{\rmdefault}{\mddefault}{\updefault}$18$}}}}}
\put(3245,350){\makebox(0,0)[lb]{\smash{{{\SetFigFont{5}{6.0}{\rmdefault}{\mddefault}{\updefault}$19$}}}}}
\put(3545,650){\makebox(0,0)[lb]{\smash{{{\SetFigFont{5}{6.0}{\rmdefault}{\mddefault}{\updefault}$20$}}}}}
\put(3545,50){\makebox(0,0)[lb]{\smash{{{\SetFigFont{5}{6.0}{\rmdefault}{\mddefault}{\updefault}$21$}}}}}
\put(3845,950){\makebox(0,0)[lb]{\smash{{{\SetFigFont{5}{6.0}{\rmdefault}{\mddefault}{\updefault}$22$}}}}}
\put(3845,350){\makebox(0,0)[lb]{\smash{{{\SetFigFont{5}{6.0}{\rmdefault}{\mddefault}{\updefault}$23$}}}}}
\put(4145,650){\makebox(0,0)[lb]{\smash{{{\SetFigFont{5}{6.0}{\rmdefault}{\mddefault}{\updefault}$24$}}}}}
\put(4520,650){\makebox(0,0)[lb]{\smash{{{\SetFigFont{5}{6.0}{\rmdefault}{\mddefault}{\updefault}$25$}}}}}
\put(4895,650){\makebox(0,0)[lb]{\smash{{{\SetFigFont{5}{6.0}{\rmdefault}{\mddefault}{\updefault}$26$}}}}}
\put(5270,650){\makebox(0,0)[lb]{\smash{{{\SetFigFont{5}{6.0}{\rmdefault}{\mddefault}{\updefault}$27$}}}}}

\put(1010,673){\ellipse{70}{70}}
\put(1010,673){\whiten\ellipse{70}{70}}
\put(1010,673){\blacken\ellipse{30}{30}}
\put(1010,673){\ellipse{30}{30}}

\put(1949,826){\ellipse{70}{70}}
\put(1949,826){\whiten\ellipse{70}{70}}
\put(1949,826){\blacken\ellipse{30}{30}}
\put(1949,826){\ellipse{30}{30}}

\put(2249,526){\ellipse{70}{70}}
\put(2249,526){\whiten\ellipse{70}{70}}
\put(2249,526){\blacken\ellipse{30}{30}}
\put(2249,526){\ellipse{30}{30}}

\put(3151,521){\ellipse{70}{70}}
\put(3151,521){\whiten\ellipse{70}{70}}
\put(3151,521){\blacken\ellipse{30}{30}}
\put(3151,521){\ellipse{30}{30}}

\put(3451,821){\ellipse{70}{70}}
\put(3451,821){\whiten\ellipse{70}{70}}
\put(3451,821){\blacken\ellipse{30}{30}}
\put(3451,821){\ellipse{30}{30}}

\put(4384,673){\ellipse{70}{70}}
\put(4384,673){\whiten\ellipse{70}{70}}
\put(4384,673){\blacken\ellipse{30}{30}}
\put(4384,673){\ellipse{30}{30}}

\put(570,673){\whiten\ellipse{50}{50}}
\put(570,673){\ellipse{50}{50}}
\put(710,673){\whiten\ellipse{50}{50}}
\put(710,673){\ellipse{50}{50}}

\put(2189,1063){\whiten\ellipse{50}{50}}
\put(2189,1063){\ellipse{50}{50}}
\put(2295,1169){\whiten\ellipse{50}{50}}
\put(2295,1169){\ellipse{50}{50}}

\put(2489,763){\whiten\ellipse{50}{50}}
\put(2489,763){\ellipse{50}{50}}
\put(2595,869){\whiten\ellipse{50}{50}}
\put(2595,869){\ellipse{50}{50}}

\put(2789,463){\whiten\ellipse{50}{50}}
\put(2789,463){\ellipse{50}{50}}
\put(2895,569){\whiten\ellipse{50}{50}}
\put(2895,569){\ellipse{50}{50}}

\put(3689,763){\whiten\ellipse{50}{50}}
\put(3689,763){\ellipse{50}{50}}
\put(3795,869){\whiten\ellipse{50}{50}}
\put(3795,869){\ellipse{50}{50}}
\put(3989,463){\whiten\ellipse{50}{50}}
\put(3989,463){\ellipse{50}{50}}
\put(4095,569){\whiten\ellipse{50}{50}}
\put(4095,569){\ellipse{50}{50}}

\put(1648,523){\whiten\ellipse{50}{50}}
\put(1648,523){\ellipse{50}{50}}
\put(1948,223){\whiten\ellipse{50}{50}}
\put(1948,223){\ellipse{50}{50}}
\put(1348,823){\whiten\ellipse{50}{50}}
\put(1348,823){\ellipse{50}{50}}

\put(3443,218){\whiten\ellipse{50}{50}}
\put(3443,218){\ellipse{50}{50}}
\put(3743,518){\whiten\ellipse{50}{50}}
\put(3743,518){\ellipse{50}{50}}
\put(4043,818){\whiten\ellipse{50}{50}}
\put(4043,818){\ellipse{50}{50}}

\end{picture}
}
}
\end{pspicture}
$$
If $x_j$ fixes the block $\om_i$, then the action on the six points of $\om_i$ is a product of three
disjoint $2$-cycles (although the generator $x_6$ in fact fixes all 162 points of $\om_3$).
Partially order the elements of
$\SS_6$ by $\sigma<\tau$ if and only if $\sigma(i)<\tau(i)$, where $i$ is the
smallest element of $\{1,\ldots,6\}$ on which $\sigma$ and $\tau$ differ.
Restrict this order to the permutations of type three disjoint $2$-cycles. The $x_j$ action on $\om_i$ is then
as,
$$
\begin{tabular}{l|l|l}
\hline
${\scriptstyle \1:1(3)2(4)3(8)4(1)7(1)}$&${\scriptstyle \1\0:2(3)4(8)5(1)}$&${\scriptstyle \1\9:1(1)2(8)4(3)}$\\
${\scriptstyle \2:1(3)2(4)3(8)7(1)}$&${\scriptstyle \1\1:2(1)3(8)7(4)}$&${\scriptstyle \2\0:1(1)2(8)}$\\
${\scriptstyle \3:1(3)2(4)5(1)7(1)}$&${\scriptstyle \1\2:5(1)7(4)}$&${\scriptstyle \2\1:1(1)2(8)3(4)4(3)}$\\
${\scriptstyle \4:1(3)4(8)5(1)}$&${\scriptstyle \1\3:2(1)3(8)4(1)7(4)}$&${\scriptstyle \2\2:1(1)2(8)3(1)5(3)}$\\
${\scriptstyle \5:1(3)3(1)4(8)5(1)}$&${\scriptstyle \1\4:3(8)7(4)}$&${\scriptstyle \2\3:1(1)2(8)3(4)}$\\
${\scriptstyle \6:3(4)4(8)5(1)}$&${\scriptstyle \1\5:1(1)3(3)5(1)7(4)}$&${\scriptstyle \2\4:1(1)2(8)5(3)}$\\
${\scriptstyle \7:4(8)5(1)}$&${\scriptstyle \1\6:3(8)4(1)7(4)}$&${\scriptstyle \2\5:1(1)4(4)5(3)7(1)}$\\
${\scriptstyle \8:2(3)3(4)4(8)5(1)}$&${\scriptstyle \1\7:1(1)7(4)}$&${\scriptstyle \2\6:3(8)4(4)5(3)7(1)}$\\
${\scriptstyle \9:2(1)5(1)7(4)}$&${\scriptstyle \1\8:1(1)4(1)7(4)}$&${\scriptstyle \2\7:2(1)3(8)4(4)5(3)7(1)}$\\
\hline
\end{tabular}
$$
using the notation ${\mathbf i}:j(k)\ldots$ to mean $x_j$ acts on the $i$-th block as the $k$-th permutation of the form
three disjoint $2$-cycles, or the identity if $j$ does not appear. The
permutations $(i,j)_k$ are all the identity in this case.
$$
\begin{pspicture}(0,0)(14,9.5)
\rput(15.125,4.15){$\star$}
\rput(.25,.5){
\rput(7,4){$\Omega_4$}
\rput(7,4.5){\BoxedEPSF{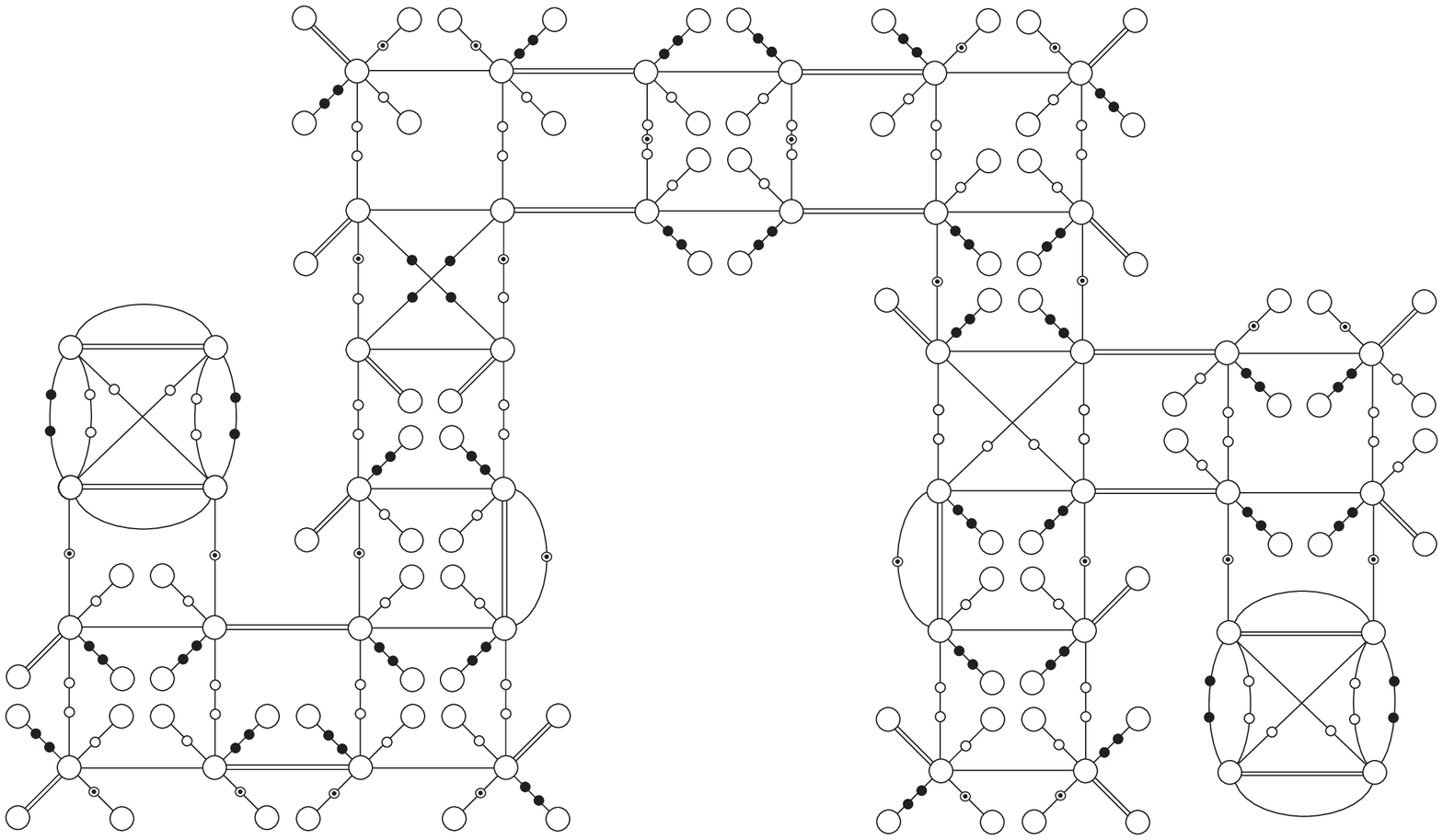 scaled 725}}
\rput(1.025,5.285){\put(0,0){\makebox(0,0)[lb]{\smash{{{\SetFigFont{5}{6.0}{\rmdefault}{\mddefault}{\updefault}$1$}}}}}}
\rput(1.02,3.655){\put(0,0){\makebox(0,0)[lb]{\smash{{{\SetFigFont{5}{6.0}{\rmdefault}{\mddefault}{\updefault}$4$}}}}}}
\rput(1.02,2.025){\put(0,0){\makebox(0,0)[lb]{\smash{{{\SetFigFont{5}{6.0}{\rmdefault}{\mddefault}{\updefault}$9$}}}}}}
\rput(.945,.395){\put(0,0){\makebox(0,0)[lb]{\smash{{{\SetFigFont{5}{6.0}{\rmdefault}{\mddefault}{\updefault}$17$}}}}}}

\rput(-.675,5.305){\put(0,0){\makebox(0,0)[lb]{\smash{{{\SetFigFont{5}{6.0}{\rmdefault}{\mddefault}{\updefault}$2$}}}}}}
\rput(-.685,3.655){\put(0,0){\makebox(0,0)[lb]{\smash{{{\SetFigFont{5}{6.0}{\rmdefault}{\mddefault}{\updefault}$5$}}}}}}
\rput(-.735,2.025){\put(0,0){\makebox(0,0)[lb]{\smash{{{\SetFigFont{5}{6.0}{\rmdefault}{\mddefault}{\updefault}$11$}}}}}}
\rput(-.755,.395){\put(0,0){\makebox(0,0)[lb]{\smash{{{\SetFigFont{5}{6.0}{\rmdefault}{\mddefault}{\updefault}$22$}}}}}}

\rput(-1.35,-.2){\put(0,0){\makebox(0,0)[lb]{\smash{{{\SetFigFont{5}{6.0}{\rmdefault}{\mddefault}{\updefault}$34$}}}}}}
\rput(-1.35,.99){\put(0,0){\makebox(0,0)[lb]{\smash{{{\SetFigFont{5}{6.0}{\rmdefault}{\mddefault}{\updefault}$35$}}}}}}
\rput(-.15,-.2){\put(0,0){\makebox(0,0)[lb]{\smash{{{\SetFigFont{5}{6.0}{\rmdefault}{\mddefault}{\updefault}$13$}}}}}}
\rput(-.15,.99){\put(0,0){\makebox(0,0)[lb]{\smash{{{\SetFigFont{5}{6.0}{\rmdefault}{\mddefault}{\updefault}$32$}}}}}}

\rput(-1.32,1.46){\put(0,0){\makebox(0,0)[lb]{\smash{{{\SetFigFont{5}{6.0}{\rmdefault}{\mddefault}{\updefault}$21$}}}}}}
\rput(-.09,1.43){\put(0,0){\makebox(0,0)[lb]{\smash{{{\SetFigFont{5}{6.0}{\rmdefault}{\mddefault}{\updefault}$6$}}}}}}
\rput(-.15,2.643){\put(0,0){\makebox(0,0)[lb]{\smash{{{\SetFigFont{5}{6.0}{\rmdefault}{\mddefault}{\updefault}$20$}}}}}}

\rput(.32,2.63){\put(0,0){\makebox(0,0)[lb]{\smash{{{\SetFigFont{5}{6.0}{\rmdefault}{\mddefault}{\updefault}$19$}}}}}}
\rput(.38,1.43){\put(0,0){\makebox(0,0)[lb]{\smash{{{\SetFigFont{5}{6.0}{\rmdefault}{\mddefault}{\updefault}$3$}}}}}}

\rput(.34,.99){\put(0,0){\makebox(0,0)[lb]{\smash{{{\SetFigFont{5}{6.0}{\rmdefault}{\mddefault}{\updefault}$30$}}}}}}
\rput(1.63,-.2){\put(0,0){\makebox(0,0)[lb]{\smash{{{\SetFigFont{5}{6.0}{\rmdefault}{\mddefault}{\updefault}$8$}}}}}}
\rput(1.56,.99){\put(0,0){\makebox(0,0)[lb]{\smash{{{\SetFigFont{5}{6.0}{\rmdefault}{\mddefault}{\updefault}$29$}}}}}}

\rput(2.055,-.2){\put(0,0){\makebox(0,0)[lb]{\smash{{{\SetFigFont{5}{6.0}{\rmdefault}{\mddefault}{\updefault}$31$}}}}}}
\rput(2.05,.99){\put(0,0){\makebox(0,0)[lb]{\smash{{{\SetFigFont{5}{6.0}{\rmdefault}{\mddefault}{\updefault}$30$}}}}}}
\rput(3.26,.99){\put(0,0){\makebox(0,0)[lb]{\smash{{{\SetFigFont{5}{6.0}{\rmdefault}{\mddefault}{\updefault}$29$}}}}}}

\rput(3.735,-.2){\put(0,0){\makebox(0,0)[lb]{\smash{{{\SetFigFont{5}{6.0}{\rmdefault}{\mddefault}{\updefault}$36$}}}}}}
\rput(3.735,.99){\put(0,0){\makebox(0,0)[lb]{\smash{{{\SetFigFont{5}{6.0}{\rmdefault}{\mddefault}{\updefault}$35$}}}}}}
\rput(4.965,.99){\put(0,0){\makebox(0,0)[lb]{\smash{{{\SetFigFont{5}{6.0}{\rmdefault}{\mddefault}{\updefault}$39$}}}}}}
\rput(4.955,-.2){\put(0,0){\makebox(0,0)[lb]{\smash{{{\SetFigFont{5}{6.0}{\rmdefault}{\mddefault}{\updefault}$32$}}}}}}

\rput(3.725,1.43){\put(0,0){\makebox(0,0)[lb]{\smash{{{\SetFigFont{5}{6.0}{\rmdefault}{\mddefault}{\updefault}$15$}}}}}}
\rput(3.725,2.62){\put(0,0){\makebox(0,0)[lb]{\smash{{{\SetFigFont{5}{6.0}{\rmdefault}{\mddefault}{\updefault}$26$}}}}}}

\rput(3.255,2.62){\put(0,0){\makebox(0,0)[lb]{\smash{{{\SetFigFont{5}{6.0}{\rmdefault}{\mddefault}{\updefault}$18$}}}}}}
\rput(3.325,1.43){\put(0,0){\makebox(0,0)[lb]{\smash{{{\SetFigFont{5}{6.0}{\rmdefault}{\mddefault}{\updefault}$7$}}}}}}

\rput(3.715,3.06){\put(0,0){\makebox(0,0)[lb]{\smash{{{\SetFigFont{5}{6.0}{\rmdefault}{\mddefault}{\updefault}$36$}}}}}}
\rput(3.725,4.25){\put(0,0){\makebox(0,0)[lb]{\smash{{{\SetFigFont{5}{6.0}{\rmdefault}{\mddefault}{\updefault}$23$}}}}}}

\rput(3.255,3.06){\put(0,0){\makebox(0,0)[lb]{\smash{{{\SetFigFont{5}{6.0}{\rmdefault}{\mddefault}{\updefault}$31$}}}}}}
\rput(3.245,4.25){\put(0,0){\makebox(0,0)[lb]{\smash{{{\SetFigFont{5}{6.0}{\rmdefault}{\mddefault}{\updefault}$12$}}}}}}

\rput(3.245,4.68){\put(0,0){\makebox(0,0)[lb]{\smash{{{\SetFigFont{5}{6.0}{\rmdefault}{\mddefault}{\updefault}$22$}}}}}}
\rput(3.715,4.675){\put(0,0){\makebox(0,0)[lb]{\smash{{{\SetFigFont{5}{6.0}{\rmdefault}{\mddefault}{\updefault}$37$}}}}}}

\rput(3.240,7.935){\put(0,0){\makebox(0,0)[lb]{\smash{{{\SetFigFont{5}{6.0}{\rmdefault}{\mddefault}{\updefault}$21$}}}}}}
\rput(2.010,7.925){\put(0,0){\makebox(0,0)[lb]{\smash{{{\SetFigFont{5}{6.0}{\rmdefault}{\mddefault}{\updefault}$25$}}}}}}
\rput(3.23,9.13){\put(0,0){\makebox(0,0)[lb]{\smash{{{\SetFigFont{5}{6.0}{\rmdefault}{\mddefault}{\updefault}$28$}}}}}}
\rput(2.005,9.14){\put(0,0){\makebox(0,0)[lb]{\smash{{{\SetFigFont{5}{6.0}{\rmdefault}{\mddefault}{\updefault}$20$}}}}}}

\rput(2.010,6.280){\put(0,0){\makebox(0,0)[lb]{\smash{{{\SetFigFont{5}{6.0}{\rmdefault}{\mddefault}{\updefault}$32$}}}}}}

\rput(2.020,3.05){\put(0,0){\makebox(0,0)[lb]{\smash{{{\SetFigFont{5}{6.0}{\rmdefault}{\mddefault}{\updefault}$11$}}}}}}

\rput(3.695,9.13){\put(0,0){\makebox(0,0)[lb]{\smash{{{\SetFigFont{5}{6.0}{\rmdefault}{\mddefault}{\updefault}$37$}}}}}}
\rput(4.910,9.13){\put(0,0){\makebox(0,0)[lb]{\smash{{{\SetFigFont{5}{6.0}{\rmdefault}{\mddefault}{\updefault}$24$}}}}}}
\rput(4.905,7.93){\put(0,0){\makebox(0,0)[lb]{\smash{{{\SetFigFont{5}{6.0}{\rmdefault}{\mddefault}{\updefault}$33$}}}}}}

\rput(6.6,9.13){\put(0,0){\makebox(0,0)[lb]{\smash{{{\SetFigFont{5}{6.0}{\rmdefault}{\mddefault}{\updefault}$13$}}}}}}
\rput(6.6,7.93){\put(0,0){\makebox(0,0)[lb]{\smash{{{\SetFigFont{5}{6.0}{\rmdefault}{\mddefault}{\updefault}$27$}}}}}}

\rput(6.6,7.5){\put(0,0){\makebox(0,0)[lb]{\smash{{{\SetFigFont{5}{6.0}{\rmdefault}{\mddefault}{\updefault}$37$}}}}}}
\rput(6.605,6.290){\put(0,0){\makebox(0,0)[lb]{\smash{{{\SetFigFont{5}{6.0}{\rmdefault}{\mddefault}{\updefault}$22$}}}}}}

\rput(7.130,9.13){\put(0,0){\makebox(0,0)[lb]{\smash{{{\SetFigFont{5}{6.0}{\rmdefault}{\mddefault}{\updefault}$8$}}}}}}
\rput(7.070,7.93){\put(0,0){\makebox(0,0)[lb]{\smash{{{\SetFigFont{5}{6.0}{\rmdefault}{\mddefault}{\updefault}$16$}}}}}}

\rput(7.0755,7.5){\put(0,0){\makebox(0,0)[lb]{\smash{{{\SetFigFont{5}{6.0}{\rmdefault}{\mddefault}{\updefault}$28$}}}}}}
\rput(7.0750,6.295){\put(0,0){\makebox(0,0)[lb]{\smash{{{\SetFigFont{5}{6.0}{\rmdefault}{\mddefault}{\updefault}$17$}}}}}}

\rput(8.755,9.125){\put(0,0){\makebox(0,0)[lb]{\smash{{{\SetFigFont{5}{6.0}{\rmdefault}{\mddefault}{\updefault}$10$}}}}}}
\rput(8.817,7.93){\put(0,0){\makebox(0,0)[lb]{\smash{{{\SetFigFont{5}{6.0}{\rmdefault}{\mddefault}{\updefault}$9$}}}}}}
\rput(10,9.13){\put(0,0){\makebox(0,0)[lb]{\smash{{{\SetFigFont{5}{6.0}{\rmdefault}{\mddefault}{\updefault}$14$}}}}}}

\rput(10.45,9.120){\put(0,0){\makebox(0,0)[lb]{\smash{{{\SetFigFont{5}{6.0}{\rmdefault}{\mddefault}{\updefault}$10$}}}}}}
\rput(10.448,7.920){\put(0,0){\makebox(0,0)[lb]{\smash{{{\SetFigFont{5}{6.0}{\rmdefault}{\mddefault}{\updefault}$11$}}}}}}
\rput(11.665,7.910){\put(0,0){\makebox(0,0)[lb]{\smash{{{\SetFigFont{5}{6.0}{\rmdefault}{\mddefault}{\updefault}$14$}}}}}}

\rput(10.465,7.485){\put(0,0){\makebox(0,0)[lb]{\smash{{{\SetFigFont{5}{6.0}{\rmdefault}{\mddefault}{\updefault}$22$}}}}}}
\rput(10.465,6.285){\put(0,0){\makebox(0,0)[lb]{\smash{{{\SetFigFont{5}{6.0}{\rmdefault}{\mddefault}{\updefault}$37$}}}}}}

\rput(11.710,6.285){\put(0,0){\makebox(0,0)[lb]{\smash{{{\SetFigFont{5}{6.0}{\rmdefault}{\mddefault}{\updefault}$38$}}}}}}

\rput(11.710,9.12){\put(0,0){\makebox(0,0)[lb]{\smash{{{\SetFigFont{5}{6.0}{\rmdefault}{\mddefault}{\updefault}$31$}}}}}}

\rput(11.80,2.595){\put(0,0){\makebox(0,0)[lb]{\smash{{{\SetFigFont{5}{6.0}{\rmdefault}{\mddefault}{\updefault}$3$}}}}}}

\rput(9.995,7.485){\put(0,0){\makebox(0,0)[lb]{\smash{{{\SetFigFont{5}{6.0}{\rmdefault}{\mddefault}{\updefault}$17$}}}}}}
\rput(10,6.285){\put(0,0){\makebox(0,0)[lb]{\smash{{{\SetFigFont{5}{6.0}{\rmdefault}{\mddefault}{\updefault}$28$}}}}}}

\rput(10.008,5.860){\put(0,0){\makebox(0,0)[lb]{\smash{{{\SetFigFont{5}{6.0}{\rmdefault}{\mddefault}{\updefault}$31$}}}}}}
\rput(8.795,5.860){\put(0,0){\makebox(0,0)[lb]{\smash{{{\SetFigFont{5}{6.0}{\rmdefault}{\mddefault}{\updefault}$14$}}}}}}
\rput(10.470,5.860){\put(0,0){\makebox(0,0)[lb]{\smash{{{\SetFigFont{5}{6.0}{\rmdefault}{\mddefault}{\updefault}$36$}}}}}}

\rput(13.40,5.855){\put(0,0){\makebox(0,0)[lb]{\smash{{{\SetFigFont{5}{6.0}{\rmdefault}{\mddefault}{\updefault}$22$}}}}}}
\rput(13.855,5.840){\put(0,0){\makebox(0,0)[lb]{\smash{{{\SetFigFont{5}{6.0}{\rmdefault}{\mddefault}{\updefault}$17$}}}}}}
\rput(15.075,5.840){\put(0,0){\makebox(0,0)[lb]{\smash{{{\SetFigFont{5}{6.0}{\rmdefault}{\mddefault}{\updefault}$10$}}}}}}

\rput(13.40,4.645){\put(0,0){\makebox(0,0)[lb]{\smash{{{\SetFigFont{5}{6.0}{\rmdefault}{\mddefault}{\updefault}$26$}}}}}}
\rput(13.85,4.645){\put(0,0){\makebox(0,0)[lb]{\smash{{{\SetFigFont{5}{6.0}{\rmdefault}{\mddefault}{\updefault}$18$}}}}}}

\rput(12.18,4.65){\put(0,0){\makebox(0,0)[lb]{\smash{{{\SetFigFont{5}{6.0}{\rmdefault}{\mddefault}{\updefault}$15$}}}}}}

\rput(12.19,4.215){\put(0,0){\makebox(0,0)[lb]{\smash{{{\SetFigFont{5}{6.0}{\rmdefault}{\mddefault}{\updefault}$14$}}}}}}
\rput(15.085,4.205){\put(0,0){\makebox(0,0)[lb]{\smash{{{\SetFigFont{5}{6.0}{\rmdefault}{\mddefault}{\updefault}$10$}}}}}}

\rput(15.15,4.63){\put(0,0){\makebox(0,0)[lb]{\smash{{{\SetFigFont{5}{6.0}{\rmdefault}{\mddefault}{\updefault}$7$}}}}}}
\rput(15.15,2.99){\put(0,0){\makebox(0,0)[lb]{\smash{{{\SetFigFont{5}{6.0}{\rmdefault}{\mddefault}{\updefault}$7$}}}}}}

\rput(13.95,2.99){\put(0,0){\makebox(0,0)[lb]{\smash{{{\SetFigFont{5}{6.0}{\rmdefault}{\mddefault}{\updefault}$9$}}}}}}
\rput(13.40,2.99){\put(0,0){\makebox(0,0)[lb]{\smash{{{\SetFigFont{5}{6.0}{\rmdefault}{\mddefault}{\updefault}$11$}}}}}}

\rput(10.49,3.02){\put(0,0){\makebox(0,0)[lb]{\smash{{{\SetFigFont{5}{6.0}{\rmdefault}{\mddefault}{\updefault}$21$}}}}}}
\rput(10.020,3.02){\put(0,0){\makebox(0,0)[lb]{\smash{{{\SetFigFont{5}{6.0}{\rmdefault}{\mddefault}{\updefault}$33$}}}}}}

\rput(10.58,2.58){\put(0,0){\makebox(0,0)[lb]{\smash{{{\SetFigFont{5}{6.0}{\rmdefault}{\mddefault}{\updefault}$8$}}}}}}
\rput(10.020,2.58){\put(0,0){\makebox(0,0)[lb]{\smash{{{\SetFigFont{5}{6.0}{\rmdefault}{\mddefault}{\updefault}$13$}}}}}}

\rput(10.495,1.37){\put(0,0){\makebox(0,0)[lb]{\smash{{{\SetFigFont{5}{6.0}{\rmdefault}{\mddefault}{\updefault}$16$}}}}}}
\rput(10.03,1.37){\put(0,0){\makebox(0,0)[lb]{\smash{{{\SetFigFont{5}{6.0}{\rmdefault}{\mddefault}{\updefault}$27$}}}}}}

\rput(10.530,-.25){\put(0,0){\makebox(0,0)[lb]{\smash{{{\SetFigFont{5}{6.0}{\rmdefault}{\mddefault}{\updefault}$20$}}}}}}
\rput(10.6,.95){\put(0,0){\makebox(0,0)[lb]{\smash{{{\SetFigFont{5}{6.0}{\rmdefault}{\mddefault}{\updefault}$3$}}}}}}
\rput(11.8,-.25){\put(0,0){\makebox(0,0)[lb]{\smash{{{\SetFigFont{5}{6.0}{\rmdefault}{\mddefault}{\updefault}$8$}}}}}}
\rput(11.74,.96){\put(0,0){\makebox(0,0)[lb]{\smash{{{\SetFigFont{5}{6.0}{\rmdefault}{\mddefault}{\updefault}$19$}}}}}}

\rput(10.020,-.25){\put(0,0){\makebox(0,0)[lb]{\smash{{{\SetFigFont{5}{6.0}{\rmdefault}{\mddefault}{\updefault}$19$}}}}}}
\rput(10.10,.95){\put(0,0){\makebox(0,0)[lb]{\smash{{{\SetFigFont{5}{6.0}{\rmdefault}{\mddefault}{\updefault}$6$}}}}}}
\rput(8.820,-.25){\put(0,0){\makebox(0,0)[lb]{\smash{{{\SetFigFont{5}{6.0}{\rmdefault}{\mddefault}{\updefault}$20$}}}}}}
\rput(8.815,.95){\put(0,0){\makebox(0,0)[lb]{\smash{{{\SetFigFont{5}{6.0}{\rmdefault}{\mddefault}{\updefault}$25$}}}}}}

\rput(2.625,8.535){\put(0,0){\makebox(0,0)[lb]{\smash{{{\SetFigFont{5}{6.0}{\rmdefault}{\mddefault}{\updefault}$31$}}}}}}
\rput(2.625,6.905){\put(0,0){\makebox(0,0)[lb]{\smash{{{\SetFigFont{5}{6.0}{\rmdefault}{\mddefault}{\updefault}$38$}}}}}}
\rput(2.625,5.275){\put(0,0){\makebox(0,0)[lb]{\smash{{{\SetFigFont{5}{6.0}{\rmdefault}{\mddefault}{\updefault}$34$}}}}}}
\rput(2.635,3.655){\put(0,0){\makebox(0,0)[lb]{\smash{{{\SetFigFont{5}{6.0}{\rmdefault}{\mddefault}{\updefault}$21$}}}}}}
\rput(2.645,2.025){\put(0,0){\makebox(0,0)[lb]{\smash{{{\SetFigFont{5}{6.0}{\rmdefault}{\mddefault}{\updefault}$16$}}}}}}
\rput(2.645,.395){\put(0,0){\makebox(0,0)[lb]{\smash{{{\SetFigFont{5}{6.0}{\rmdefault}{\mddefault}{\updefault}$28$}}}}}}

\rput(4.3,8.535){\put(0,0){\makebox(0,0)[lb]{\smash{{{\SetFigFont{5}{6.0}{\rmdefault}{\mddefault}{\updefault}$36$}}}}}}
\rput(4.3,6.910){\put(0,0){\makebox(0,0)[lb]{\smash{{{\SetFigFont{5}{6.0}{\rmdefault}{\mddefault}{\updefault}$40$}}}}}}
\rput(4.305,5.275){\put(0,0){\makebox(0,0)[lb]{\smash{{{\SetFigFont{5}{6.0}{\rmdefault}{\mddefault}{\updefault}$39$}}}}}}
\rput(4.315,3.655){\put(0,0){\makebox(0,0)[lb]{\smash{{{\SetFigFont{5}{6.0}{\rmdefault}{\mddefault}{\updefault}$33$}}}}}}
\rput(4.335,2.025){\put(0,0){\makebox(0,0)[lb]{\smash{{{\SetFigFont{5}{6.0}{\rmdefault}{\mddefault}{\updefault}$27$}}}}}}
\rput(4.345,.395){\put(0,0){\makebox(0,0)[lb]{\smash{{{\SetFigFont{5}{6.0}{\rmdefault}{\mddefault}{\updefault}$37$}}}}}}

\rput(6,8.535){\put(0,0){\makebox(0,0)[lb]{\smash{{{\SetFigFont{5}{6.0}{\rmdefault}{\mddefault}{\updefault}$26$}}}}}}
\rput(6,6.9){\put(0,0){\makebox(0,0)[lb]{\smash{{{\SetFigFont{5}{6.0}{\rmdefault}{\mddefault}{\updefault}$35$}}}}}}

\rput(7.675,8.525){\put(0,0){\makebox(0,0)[lb]{\smash{{{\SetFigFont{5}{6.0}{\rmdefault}{\mddefault}{\updefault}$18$}}}}}}
\rput(7.68,6.905){\put(0,0){\makebox(0,0)[lb]{\smash{{{\SetFigFont{5}{6.0}{\rmdefault}{\mddefault}{\updefault}$29$}}}}}}

\rput(9.375,8.525){\put(0,0){\makebox(0,0)[lb]{\smash{{{\SetFigFont{5}{6.0}{\rmdefault}{\mddefault}{\updefault}$19$}}}}}}
\rput(9.375,6.875){\put(0,0){\makebox(0,0)[lb]{\smash{{{\SetFigFont{5}{6.0}{\rmdefault}{\mddefault}{\updefault}$30$}}}}}}

\rput(11.08,8.525){\put(0,0){\makebox(0,0)[lb]{\smash{{{\SetFigFont{5}{6.0}{\rmdefault}{\mddefault}{\updefault}$20$}}}}}}
\rput(11.08,6.875){\put(0,0){\makebox(0,0)[lb]{\smash{{{\SetFigFont{5}{6.0}{\rmdefault}{\mddefault}{\updefault}$32$}}}}}}

\rput(9.405,5.26){\put(0,0){\makebox(0,0)[lb]{\smash{{{\SetFigFont{5}{6.0}{\rmdefault}{\mddefault}{\updefault}$25$}}}}}}
\rput(9.415,3.63){\put(0,0){\makebox(0,0)[lb]{\smash{{{\SetFigFont{5}{6.0}{\rmdefault}{\mddefault}{\updefault}$23$}}}}}}
\rput(9.420,2.0){\put(0,0){\makebox(0,0)[lb]{\smash{{{\SetFigFont{5}{6.0}{\rmdefault}{\mddefault}{\updefault}$15$}}}}}}
\rput(9.425,.360){\put(0,0){\makebox(0,0)[lb]{\smash{{{\SetFigFont{5}{6.0}{\rmdefault}{\mddefault}{\updefault}$14$}}}}}}

\rput(11.090,5.26){\put(0,0){\makebox(0,0)[lb]{\smash{{{\SetFigFont{5}{6.0}{\rmdefault}{\mddefault}{\updefault}$24$}}}}}}
\rput(11.095,3.63){\put(0,0){\makebox(0,0)[lb]{\smash{{{\SetFigFont{5}{6.0}{\rmdefault}{\mddefault}{\updefault}$12$}}}}}}
\rput(11.170,2.0){\put(0,0){\makebox(0,0)[lb]{\smash{{{\SetFigFont{5}{6.0}{\rmdefault}{\mddefault}{\updefault}$7$}}}}}}
\rput(11.130,.360){\put(0,0){\makebox(0,0)[lb]{\smash{{{\SetFigFont{5}{6.0}{\rmdefault}{\mddefault}{\updefault}$10$}}}}}}

\rput(12.765,5.245){\put(0,0){\makebox(0,0)[lb]{\smash{{{\SetFigFont{5}{6.0}{\rmdefault}{\mddefault}{\updefault}$13$}}}}}}
\rput(12.855,3.620){\put(0,0){\makebox(0,0)[lb]{\smash{{{\SetFigFont{5}{6.0}{\rmdefault}{\mddefault}{\updefault}$6$}}}}}}

\rput(14.525,5.235){\put(0,0){\makebox(0,0)[lb]{\smash{{{\SetFigFont{5}{6.0}{\rmdefault}{\mddefault}{\updefault}$8$}}}}}}
\rput(14.545,3.605){\put(0,0){\makebox(0,0)[lb]{\smash{{{\SetFigFont{5}{6.0}{\rmdefault}{\mddefault}{\updefault}$3$}}}}}}

\rput(12.875,1.965){\put(0,0){\makebox(0,0)[lb]{\smash{{{\SetFigFont{5}{6.0}{\rmdefault}{\mddefault}{\updefault}$2$}}}}}}
\rput(14.555,1.965){\put(0,0){\makebox(0,0)[lb]{\smash{{{\SetFigFont{5}{6.0}{\rmdefault}{\mddefault}{\updefault}$1$}}}}}}

\rput(12.875,.335){\put(0,0){\makebox(0,0)[lb]{\smash{{{\SetFigFont{5}{6.0}{\rmdefault}{\mddefault}{\updefault}$5$}}}}}}
\rput(14.555,.335){\put(0,0){\makebox(0,0)[lb]{\smash{{{\SetFigFont{5}{6.0}{\rmdefault}{\mddefault}{\updefault}$4$}}}}}}
}
\end{pspicture}
$$
Finally, diagram $\om_4$ above gives a $\Gamma$-action on $40$ blocks of imprimitivity, each
of size eight, where repeated vertices and edges are to be thought of as being identified. 
No block is fixed by any generator except for $x_6$, which, as for $\om_3$, 
fixes all $320$ points.
The permutations $(i,j)_k$ not equal to the identity are given by,
$$
\begin{tabular}{l|l|l}
\hline
${\scriptstyle (1,2)_2,(1,4)_4,(12,24)_4= (2,3)(4,6)(5,7)}$&${\scriptstyle
(2,4)_7= (2,4)(3,5)(6,7)}$&${\scriptstyle (1,5)_7,(27,16)_1,(8,13)_1= (3,6,4,7,5)}$\\

${\scriptstyle (2,5)_4,(4,5)_2=(2,4,6,3,5)}$&${\scriptstyle
(1,4)_5,(1,2)_1= (3,6,5,4)}$&${\scriptstyle (5,11)_3,(10,19)_5= (6,7)}$\\

${\scriptstyle (11,20)_7= (5,6)}$&${\scriptstyle (9,19)_7,(15,23)_2,(14,25)_2,(6,14)_7= (3,4)}$&${\scriptstyle
(11,21)_2,(11,22)_4= (2,3,4)(6,8,7)}$\\

${\scriptstyle (9,17)_4= (2,3)(5,8,7,6)}$&${\scriptstyle
(22,35)_5,(21,33)_1= (2,3,4)(5,6,7)}$&${\scriptstyle
(22,34)_2,(21,34)_4= (3,4,5)(6,7)}$\\

${\scriptstyle (22,17)_1,(10,3)_7= (3,4)(5,7,6)}$&${\scriptstyle (17,30)_7= (5,6,7,8)}$&${\scriptstyle 
(2,5)_5,(4,5)_1,(37,28)_1= (2,3,4)(5,7,6)}$\\

${\scriptstyle
(4,9)_3,(38,31)_4= (6,7,8)}$&${\scriptstyle (9,11)_1,(30,32)_1= (2,3,4)(5,6,8,7)
}$&${\scriptstyle (22,32)_7,(20,19)_1,(23,15)_3= (2,3)(6,7,8)
}$\\

${\scriptstyle (17,28)_2,(10,20)_3= (4,5)(6,7,8)}$&${\scriptstyle (28,30)_5= (1,3)(5,6,7)}$&
${\scriptstyle (28,31)_3= (1,2)(3,4)(5,6)}$\\

${\scriptstyle (28,29)_7= (1,2,3)(4,6,5)(7,8)}$&${\scriptstyle (16,28)_4,(29,17)_5,(23,12)_1= (2,3)(4,5,6)}$&${\scriptstyle
(27,37)_4,(26,36)_2= (4,5)}$ \\

${\scriptstyle (16,18)_7= (1,2)(3,5,4)}$&${\scriptstyle (16,21)_3= (3,5)(4,6,7)}$&${\scriptstyle
(27,33)_2,(26,35)_4= (3,5)(4,6)(7,8)}$\\

${\scriptstyle (27,33)_3,(26,35)_3= (2,3,5,6,7,4)}$&${\scriptstyle
(21,31)_7= (2,4,3)(7,8)}$&${\scriptstyle (33,36)_7,(35,37)_7= (1,2)(3,4)(5,6,7,8)}$\\

${\scriptstyle (33,39)_4,(35,40)_2= (3,4,6,5)(7,8)}$&${\scriptstyle
(34,39)_1,(34,40)_5= (2,4,6,3,5)(7,8)}$&${\scriptstyle
(38,39)_5,(38,40)_1= (1,2,3)(5,7,6)}$\\

${\scriptstyle (34,38)_7= (2,4)(3,5)(6,7,8)}$&${\scriptstyle (34,38)_3= (1,2,4,5,6,3)(7,8)}$&${\scriptstyle
(39,40)_7= (1,2)(3,6)(5,7)}$\\

${\scriptstyle (39,40)_3= (1,3)(2,5)(6,7)}$&${\scriptstyle
(31,36)_1= (1,2,3)(5,7)(6,8)}$&${\scriptstyle (36,37)_3= (1,3)(2,4)(6,8)}$\\

${\scriptstyle  (2,6)_3= (7,8)}$&${\scriptstyle
(26,27)_7= (1,2)(3,6)(4,7)}$&${\scriptstyle (29,35)_1,(13,24)_2= (4,5,6)}$\\

${\scriptstyle (18,29)_3= (2,3,4)(6,7)}$&${\scriptstyle
(18,29)_4= (3,4,5)(7,8)}$&${\scriptstyle (18,19)_2,(8,10)_2= (1,2)(4,5)}$\\

${\scriptstyle (29,30)_2= (1,3,2)(7,8)}$&${\scriptstyle (20,31)_2,(6,3)_1= (5,7,6)}$&${\scriptstyle (19,30)_4= (2,3,4))}$\\

${\scriptstyle (20,32)_4= (3,4)(5,6)}$&${\scriptstyle (32,38)_2= (3,4)(7,8)}$&${\scriptstyle (32,37)_5=(1,2,3,4)(5,7,8,6)}$\\

${\scriptstyle (25,30)_3= (1,2,4,3)(6,7,8)}$&${\scriptstyle (24,32)_3= (1,2,3)(6,7)}$&${\scriptstyle (25,31)_5= (2,3,4)(5,7,8,6)}$\\

${\scriptstyle
(24,36)_5= (3,4)(5,7)}$&${\scriptstyle (24,25)_1= (1,2)(5,6,7)}$&${\scriptstyle (23,25)_4= (1,2)(4,5)(6,7)}$\\

${\scriptstyle
(23,24)_7= (1,2)(4,5,6)}$&${\scriptstyle (12,25)_7,(12,6)_2= (2,3)(5,6,7)}$&${\scriptstyle (23,33)_5= (4,6)(7,8)}$\\

${\scriptstyle
(12,21)_5= (3,4)(6,7,8)}$&${\scriptstyle (7,12)_3= (4,5,6,8,7)}$&${\scriptstyle (7,8)_7= (1,2)(3,6)(4,5,7)}$\\

${\scriptstyle
(15,27)_5= (3,4,5)}$&${\scriptstyle (7,16)_5= (3,7,5,6,4)}$&${\scriptstyle (14,15)_4,(10,7)_4= (1,2)(3,5)(6,7)}$\\

${\scriptstyle
(14,20)_5= (2,3)(7,8)}$&${\scriptstyle (14,19)_3= (2,3)(4,5)(6,8,7)}$&${\scriptstyle (13,22)_3= (3,5,6,7,4)}$\\

${\scriptstyle
(8,17)_3,(7,3)_2= (3,5,6,4)}$&${\scriptstyle (13,15)_7= (1,2)(3,4)}$&${\scriptstyle (13,26)_5= (3,5,6)(4,7)}$\\

${\scriptstyle
(6,13)_4= (4,5)(6,7)}$&${\scriptstyle (3,8)_4= (3,5,4,7,6)}$&${\scriptstyle (6,11)_5= (2,3,4)(5,6)(7,8)}$\\

${\scriptstyle
(3,9)_5= (5,7)}$&${\scriptstyle (1,3)_3= (3,5,4)(6,7,8)}$&\\
\hline
\end{tabular}
$$
As per usual, let $U=\otimes U_i$, having checked that the $U_i$ collectively avoid
the torsion in $\Gamma$, so that we obtain a (non-compact) $6$-dimensional manifold $M^6$.

To compute the volume of $M^6$, let $\Omega$ be a $\Gamma$-orbit on
$\om_2\times \om_3$ arising from $U_2\otimes U_3$. By Lemma
\ref{lemma3.2} we have that $2^6\,3^4$ divides $|\om|$, while on the
other hand, 
$|\om|\leq 2^7\,3^4$. The element 
$x_1$
of order $2$ fixes no point of $\om_3$ but fixes the point $\dag$ of
$\om_2$, thus the subgroup $F$ generated by $x_1$ does so. By
Proposition \ref{prop2}, we have that $2|\om_2|=2^7$ divides $|\om|$, giving
$|\om|=2^7\,3^4=|\om_2|\times |\om_3|$. Thus $U_2\otimes U_3$ is a
transitive $\C\Gamma$-module.

Let $\ov{\om}_4$ be the blocks of imprimitivity of the $\Gamma$-action on
$\om_4$. The element $\gamma=x_1x_3x_2x_4$ of order $5$ fixes the point $(\dag,*)$ of
$\om_2\times\om_3$ but acts on $\ov{\om}_4$ as a product of eight
$5$-cycles ($*$ is a point of $\om_3$ in the first block of imprimitivty).
Thus $F=\langle\gamma\rangle$ satisfies the conditions of
Proposition \ref{prop2}, hence the $\Gamma$-orbits on
$\om_2\times\om_3\times\ov{\om}_4$ have size divisible by $5(|\om_2|\times
|\om_3|)=2^7\,3^4\,5$. The elements,
$$
\gamma=(x_5x_4)^2x_5x_1,\delta=x_3x_2x_1x_5x_7x_3x_5x_3x_4x_1x_2x_3
$$
fix the point $(\dag,*)$ of $\om_2\times\om_3$ (where $*$ is in the
first block of imprimitivity). The faithful action of
$\gamma$ and $\delta$ on the roots for $E_6$ are as permutations
satisfying $\gamma^2=\delta^2=(\gamma\delta)^4=1$ so that
$F=\langle\gamma,\delta\rangle$ hass order $\geq 8$. Moreover,
$\gamma$ and $\delta$ act as a product of twenty $2$-cycles and
$\gamma\delta$ as a product of ten $4$-cycles on $\ov{\om}_4$, giving that
$F$ satisfies the conditions of Proposition \ref{prop2}. Thus
$2^3(|\om_2|\times|\om_3|)$ divides the size of any $\Gamma$-orbit on
$\om_2\times\om_3\times\ov{\om}_4$. Hence $2^{10}\,3^4\,5$ divides this
size which in turn is $\leq 2^{10}\,3^4\,5$, and so we have equality. Thus
$U_2\otimes U_3\otimes \ov{U}_4$ is transitive. 

Similarly $x_1x_2x_3x_4x_7$ has order $8$, acting on $\om_4$ as fourty
$8$-cycles while fixing the point $(\dag,*,\star)$ of 
$\om_2\times\om_3\times\ov{\om}_4$, hence $U_2\otimes U_3\otimes U_4$ is
transitive ($*$ as before is in the first block of imprimitivity of $\om_3$). 
Finally, 
$$
\mu=x_3x_4x_5x_6x_5x_4x_3x_7x_3x_4x_5x_6x_5x_4x_3x_4x_5x_6x_5x_4x_5x_6x_5x_6x_7,
$$
is an element of finite order, being in the copy of $B_5$ generated by 
$x_3,\ldots, x_7$.
As it involves an odd number of generators, it acts as a $2$-cycle on $\om_1$ so
in particular
is non-trivial
(actually it turns out to be an involution). Thus $\langle \mu\rangle$ has
finite order $\geq 2$. As $x_6$ fixes
all the points in both $\om_3$ and $\om_4$, the action of $\mu$ on these sets is the same as that of 
the element obtained by removing the occurances of $x_6$ from it. As the resulting word collapses to the 
identity  we have that $\mu$ fixes every point of $\om_3$ and $\om_4$. 
As it fixes the point $\dag$ of $\om_2$
also, we have $\langle \mu\rangle$ fixing a point in $\om_2\times\om_3\times\om_4$.
Proposition \ref{prop2}
gives $2|\om_2\times\om_3\times\om_4|=\times |\om_i|$ divides
the size of a $\Gamma$-orbit, and so $U=\otimes U_i$ is transitive.

Thus $M^6$ has Euler characteristic
$\chi(\Gamma)\prod|\om_i|=(-1/\LL(\Gamma))\times 16\LL(\Gamma)=-16$. As $U_1$ is clearly an orientable module,
$U$ is as well by Proposition \ref{prop3}, so that $M^6$ is orientable.



{}

\end{document}